\documentclass[review,10pt]{elsarticle}
\usepackage[utf8]{inputenc}
\usepackage{microtype}
\usepackage[margin=2.5cm]{geometry}
\usepackage{setspace}
\usepackage[none]{hyphenat}
\usepackage[pagewise]{lineno}
\usepackage[english]{babel}
\usepackage{color}
\usepackage{booktabs}
\usepackage{longtable}
\usepackage[most]{tcolorbox}
\usepackage{ragged2e}
\usepackage[overload]{empheq}
\usepackage[justification=justified]{caption}
\usepackage{enumitem}

\usepackage{amsmath,amsfonts,amssymb,mathtools}
\usepackage{cases}
\usepackage{multirow}
\usepackage{graphicx,float}
\usepackage{algpseudocode}
\usepackage{subcaption}
\usepackage{eurosym}
\usepackage{amstext} 
\usepackage{arydshln}
\usepackage{flushend}
\usepackage{threeparttable}
\usepackage{url}
\usepackage{colortbl}
\usepackage{tikz}
\usetikzlibrary{arrows.meta, positioning, calc}
\usepackage{tabularx}
\usepackage{framed}
\newcolumntype{L}{>{\raggedright\arraybackslash}X}
\usepackage{soul}
\usepackage{etoolbox}
\usepackage{xspace}

\usepackage[]{hyperref}
\hypersetup{colorlinks, linkcolor=blue, citecolor=blue, urlcolor=blue}

\usepackage[ruled,linesnumbered]{algorithm2e}
\DeclareRobustCommand{\officialeuro}{%
  \ifmmode\expandafter\text\fi
  {\fontencoding{U}\fontfamily{eurosym}\selectfont e}}
\usepackage{framed} 
\immediate\write18{%
makeindex Main_File.nlo -s nomencl.ist -o Main_File.nls 
}
\usepackage{multicol} 
\usepackage{nomencl}[noprefix] 
\makenomenclature
\renewcommand*\nompreamble{\begin{multicols}{2}}
\renewcommand*\nompostamble{\end{multicols}}

\def\tsc#1{\csdef{#1}{\textsc{\lowercase{#1}}\xspace}}
\tsc{WGM}
\tsc{QE}
\tsc{EP}
\tsc{PMS}
\tsc{BEC}
\tsc{DE}

\begin{document}

    \justifying
    \begin{frontmatter}   
        \title{A Bilevel Optimization Model for Bottom-Up Coordination of Multiple Low-Voltage Energy Communities and the Medium-Voltage Network}

        \author[1,2]{Fernando García-Muñoz\corref{mycorrespondingauthor}}
        \ead{fernando.garciam@usach.cl}
        \cortext[mycorrespondingauthor]{Corresponding author}

        \author[1,2]{Sebastián Dávila}
        \author[1,2]{Luis Rojo-González}
        \author[1,2]{Cristian Duran-Mateluna}

        \address[1]{University of Santiago of Chile (USACH), Faculty of Engineering, Industrial Engineering Department, Chile}
        \address[2]{University of Santiago of Chile (USACH), Faculty of Engineering, Program for the Development of Sustainable Production Systems (PDSPS), Chile}

    \begin{abstract}
        The increasing penetration of distributed energy resources (DERs) is transforming low-voltage (LV) networks into active systems, including energy communities, whose generation, storage, and energy exchange activities require enhanced coordination with upstream medium-voltage (MV) networks. In the proposed Stackelberg structure, a local market operator aggregates multiple LV communities and acts as a single leader, determining DER operations and boundary energy exchanges, while the MV network serves as the follower, ensuring efficient system feasibility through an economic dispatch that includes both conventional and utility-scale PV generation. The proposed bottom-up coordination scheme concentrates discrete DER scheduling at the LV level while the MV level retains a convex continuous formulation, enabling an exact single-level reformulation via the Karush–Kuhn–Tucker (KKT) conditions. In addition, a distributed coordination algorithm that combines Lagrangian Dual Decomposition (LDD) with the Alternating Direction Method of Multipliers (ADMM) is developed to coordinate LV communities in parallel while preserving data confidentiality. The framework is validated using the IEEE 33-bus system at the MV level and six European 206-bus LV test feeders. Results indicate that the LDD-ADMM algorithm closely matches the exact reformulation, with an average relative deviation of $1.7 \times 10^{-4}$, with deviations confined to periods of scarcity for the cheap resource. Furthermore, leaders' decisions can induce operating conditions that increase followers' costs relative to their independently optimal dispatch, a pattern reinforced by comparison with a feasibility-based single-level relaxation that satisfies the required energy exchanges but fails to achieve a cost-efficient allocation of MV resources.

        

    \end{abstract}
    
    \begin{keyword}
    \emph{Bilevel optimization; ADMM-based decomposition; Distributed energy resources; Multi-voltage coordination}.
    \end{keyword}

\end{frontmatter}

\section{Introduction}\label{sec:Intro}
The increasing deployment of distributed energy resources (DERs), such as rooftop photovoltaic (PV) systems and battery energy storage systems (BESS), is reshaping the operational structure of distribution networks (DNs) \cite{ADHAM20251980}. Particularly at the Low-Voltage (LV) level, the proliferation of DERs is driving a shift from passive end-consumers to active prosumers who, through emerging local or flexibility markets, could participate in power system operations, increasing the complexity of modeling DN behavior \cite{9714405}. As these challenges intensify, coordination mechanisms are required to ensure secure and efficient interoperability between these new actors and resources across different voltage levels, properly managing local decision-making and system-wide operational constraints \cite{en14154451}.

In this context, numerous approaches have been proposed in the literature to enhance interaction across different voltage levels and enable efficient utilization of DERs in DNs \cite{en15197312}. Thus, most of these studies adopt a coordination scheme in which the Transmission system operator (TSO) remains the primary decision-making authority, supervising downstream operations and relying on aggregated information from LV levels. Within this paradigm, several studies have explored different operational functions and services that could be procured from these new active distribution networks (ADN), including unit commitment~\cite{9006872}, flexibility management~\cite{10609023,9939101}, frequency regulation~\cite{1664986}, voltage control~\cite{9543104}, and market integration of distributed storage~\cite{4956966}. In parallel, practical initiatives have complemented these academic developments, such as SmartNet \cite{migliavacca2017smartnet} and CoordiNet projects \cite{GARCIAMUNOZ2023109386}, which have provided empirical insights by testing different coordination schemes and market-based mechanisms for flexibility procurement across MV and LV grids. Similarly, in a Swiss case study presented in~\cite{KALANTARNEYESTANAKI2024110747}, DSOs deliver local flexibility in response to TSO-defined power needs under centralized coordination, while the EUniversal project \cite{mourao2023euniversal} has furthered this trend by proposing a universal interface to enable more streamlined, decentralized flexibility activation.

However, as LV networks become more dynamic and populated with a heterogeneous set of DERs, such hierarchical schemes, where the TSO leads the coordination, may struggle to adequately capture the increasing complexity of these systems, potentially leading to inefficient dispatch, underutilization of local flexibility, or even technical issues when fast local dynamics are not properly accounted for \cite{UZUM2024109976}. Specifically, from an optimization modeling perspective, this complexity arises from two main sources: (i) the hierarchical coupling across voltage levels, where upper layers (typically the TSO) impose operational setpoints or constraints that must be complied with by lower layers (DSOs and end-users), and (ii) the increasing presence of discrete decision variables associated with BESS operation, electric vehicles (EV) charging, peer-to-peer (P2P) interactions, and flexibility service activation. Under a traditional \textit{top-down} structure in which the TSO establishes system-wide setpoints first, most of the discrete decision-making is pushed to the DSOs and end-users at lower levels, making the resulting problems inherently difficult to solve with conventional duality-based approaches, such as the Karush–Kuhn–Tucker (KKT) reformulation or decomposition algorithms \cite{math14040631}. 

An alternative coordination scheme recently explored in the literature inverts this decision-making order by allowing end-users or LV networks to act first, determining their local strategies based on individual capabilities and constraints, while the upstream operator aggregates these decisions and ensures global feasibility \cite{ZHANG2024123073, 10943263}. This \textit{bottom-up} structure distributes complexity among multiple decision-makers and enables explicit consideration of DER from the outset of the decision-making process. By shifting complexity to the first movers, the coordination problem could become more tractable, enabling the application of traditional optimization methods and better aligning with the decentralized, heterogeneous nature of modern distribution grids.

Although these studies introduce an alternative bottom-up coordination scheme, they remain constrained in scope. For example, these last contributions rely on a two-stage hierarchical sequencing, which results in independent optimization problems coordinated only through coupling constraints, rather than explicitly modeling the hierarchical interaction as a bilevel structure that captures the leader–follower nature and its strategic interdependencies. Furthermore, these closely related studies do not model interactions between multiple LV areas that act as energy communities and their upstream MV network, nor do they consider intercommunity energy exchanges. In this regard, this paper proposes a bilevel optimization model to represent a bottom-up coordination scheme in which a market operator aggregates multiple LV energy communities and acts as the single leader in a Stackelberg game, while the MV network operator responds as the follower to ensure system feasibility. To the best of the authors’ knowledge, this is among the first formulations to capture the coordinated operation of multiple LV communities with an MV network within a bilevel bottom-up framework. The resulting bilevel problem is solved via a KKT-based reformulation of the follower problem and a distributed coordination algorithm that combines Lagrangian relaxation and the Alternating Direction Method of Multipliers (ADMM), enabling the parallel resolution of LV-level decisions while handling the privacy of individual communities.

The remainder of the paper is organized as follows. Section \ref{sec:Literature_review} reviews the related literature on coordination across voltage levels and details the main contributions of this work. Section \ref{sec:Coordination_scheme} introduces the proposed bottom-up coordination scheme and the corresponding bilevel optimization model. Section \ref{sec:Resolution_strategy} presents the decomposition coordination algorithm. Section \ref{sec:Case_study} describes the case study and discusses the computational results. Section \ref{sec:Conclusions} concludes the paper.

\section{Literature Review}\label{sec:Literature_review}
The literature on coordination across voltage levels has grown substantially in recent years, with most existing studies focusing on the TSO–DSO interface \cite{Yang_review_2026}. However, in the absence of a unified accepted taxonomy to systematically classify these approaches, this work organizes the existing literature around three features commonly found in current contributions: the \textit{hierarchical structure} of the decision-making problem, the \textit{coordination scheme} between system operators, and the \textit{solution method} used for its implementation, thereby employing a unified categorization to facilitate comparative analysis across studies. 

The \textit{hierarchical structure} reflects the sequential nature of decision-making due to differences in voltage levels across power systems, and is typically represented either as a two-stage process with independent optimization problems or as a bilevel formulation with explicit leader–follower coupling. The \textit{coordination scheme} is defined by the sequencing and directional structure of decision-making between operators. Accordingly, different interaction patterns arise depending on which entity initiates the decision process and how many agents are involved. In the TSO–DSO interface, three schemes are commonly reported in the literature: (i) TSO-led coordination, where the transmission operator determines activation decisions and the distribution operator subsequently ensures local feasibility; (ii) DSO-led coordination, where the distribution operator allocates flexibility first and the transmission operator relies on the resulting aggregated availability; and (iii) joint sequential coordination, where both operators iteratively exchange decisions and validations to reach a feasible coordinated outcome. This classification can be naturally extended to more general settings, including interactions within distribution systems, where additional actors such as aggregators or flexibility providers lead to one-to-many or many-to-one decision relations, as captured in the \textit{Decision Sequence} column in Table \ref{Literature_Review_Table}. Finally, the \textit{solution method} specifies how the coordination scheme is implemented in practice: a centralized approach solves a single optimization problem with full information sharing, whereas a distributed approach decomposes the problem into local subproblems and coordinates them through iterative information exchange.

Following the above categorization, early works, such as \cite{YUAN2017600}, initiated TSO-led coordination through a two-stage hierarchy in which the transmission operator first sets the dispatch and aggregated power requirements for DSOs. However, DERs are not explicitly modeled, and a convex second-order cone programming (SOCP) relaxation is used for the transmission network (TN), which is uncommon given the meshed configuration of the transmission level. Along the same lines, authors in \cite{9503337} propose a TSO-led two-stage coordination approach, solved through a surrogate Lagrangian relaxation, that iteratively coordinates independent TSO and DSO subproblems. Yet, DERs are aggregated only at the DSO level without detailed modeling, and the approach primarily focuses on improving computational efficiency rather than capturing distribution-level complexity. Other works incorporate the complexity of DNs with high DER penetration, including PV, BESS, and EVs. Specifically, the authors in \cite{9615006} propose a TSO-led scheme in which the DSO response is learned as a price-sensitive aggregate bid curve to represent flexible demand, thereby avoiding explicit DN modeling. Likewise, the works in \cite{10202840} and \cite{10124221} adopt ADMM to solve distributed TSO-led coordination schemes that explicitly include DERs at the DN level. However, these formulations relax the binary charging/discharging decisions of BESS and EVs into continuous variables to preserve problem tractability under the ADMM framework, which may result in simultaneous charging and discharging at certain hours. Additionally, the authors in \cite{10892311} extend the TSO-led approach under a bilevel coordination scheme, incorporating AC network constraints and a single DSO. They apply the Analytical Target Cascading (ATC) method to address coordination issues when integrating PV and BESS at the distribution level. However, as in previous studies, BESS and EV operations are modeled as continuous variables. In addition, only in some cases are emerging local market structures or flexibility markets considered, further limiting the representation of their interaction within DNs. This limitation becomes increasingly relevant as the penetration of these resources grows, since their proper representation in the coordination scheme requires additional modeling strategies to ensure realistic operational behavior.

In this context, some studies have addressed these modeling challenges by focusing on the internal coordination mechanisms within the DN, where multiple actors, such as DSOs, aggregators, local market operators, and end-users across MV and LV levels, interact and exchange energy or flexibility. To capture these relationships, hierarchical structures such as two-stage or bilevel coordination schemes are also adopted, enabling the explicit representation of bidding, flexibility provision, and operational responses among the different entities within the DN. For example, in \cite{bahramara2015modelling} and \cite{BAHRAMARA2016169}, the authors propose a hierarchical bilevel optimization framework to model operational coordination among the distribution company (the leader) and multiple microgrids MGs (the followers). In this formulation, the DSO maximizes its profit while each MG minimizes its operational cost, establishing a local retail electricity market structure. The resulting nonlinear bilevel models are reformulated as single-level optimization problems using KKT conditions and dual theory, enabling their solution as mixed-integer linear programs. In the same bilevel optimization line, the authors in \cite{FERNANDEZ2021123254} propose a community energy management system (EMS) that coordinates P2P energy trading among prosumers interconnected with a central BESS and the utility grid. Their approach formulates the interaction as a Stackelberg game, with the BESS as the leader and consumers as followers, and reformulates the bilevel problem into a single-level problem via KKT conditions. However, these three initial bilevel approaches \cite{bahramara2015modelling,BAHRAMARA2016169,FERNANDEZ2021123254} do not incorporate network operational constraints in their modeling frameworks.

In contrast, DN limitations have been considered in hierarchical coordination approaches where the distribution level takes a more central role. One of the earliest contributions in this direction is the work of \cite{LECADRE2019317}, who analyzes different TSO–DSO coordination schemes, including a Stackelberg formulation in which the DSO acts as the leader and explicitly incorporates both TN and DN constraints. However, their study is limited to a single DSO and a small-scale DN, and does not consider BESS integration or emerging local energy markets. Similarly, the authors in \cite{GAZAFROUDI2021117575} propose a two-stage hierarchical coordination framework for energy and flexibility trading within a single DN, in which a local market operator coordinates P2P transactions in the first stage and interacts with the DSO in the second stage to dispatch flexibility and maintain network constraints. Likewise, the work in \cite{IRIA2022122962} proposes a two-stage coordination framework for MV–LV DNs, in which an aggregator submits bids to real-time energy and reserve markets in the first stage, while the DSO validates network feasibility in the second stage through iterative ADMM-based negotiation. However, the approach remains focused on a single aggregator and a single DSO within a single DN. Overall, while these works highlight a clear shift toward hierarchical coordination at the distribution level, existing models either neglect detailed network constraints or remain limited to small-scale.

In recent years, bilevel approaches have been tested in slightly larger case studies that incorporate multiple followers and, in some cases, DNs limitations. For example, \cite{9591485} develops a bilevel coordination scheme between the DSO and multiple EV charging stations, where the DSO imposes network constraints at the upper level while each station optimizes EV charging and P2P energy exchange at the lower level. Although EV flexibility and P2P interactions are represented, binary charging decisions are relaxed in the lower level to enable an ADMM-based solution. Along the same line, \cite{9645337} examines P2P energy trading among multiple prosumers through a bilevel model in which DN operational constraints are enforced at the upper level, while lower-level agents independently schedule their DERs and trading. A similar relaxation strategy is also adopted in \cite{9721228}, where a central market operator coordinates energy exchanges across local markets, although DN network limitations are not considered in their formulation. Additionally, the authors in \cite{WANG2023109065} incorporate uncertainty into a bilevel P2P coordination scheme among multiple MGs, which is solved using a hybrid ATC–ADMM approach. However, each MG is still modeled individually, limiting the representation of collective or discrete flexibility behaviors. The authors in \cite{JIA2024123259} also adopt an ADMM-based distributed solution in which a market operator acts as leader coordinating energy scheduling, while buildings operate individually as followers optimizing their internal resources under P2P-like interactions. Overall, while these approaches capture emerging distribution-level complexity, they rely on continuous relaxations of discrete flexibility decisions, which may yield operationally unrealistic solutions, highlighting the need for more advanced hierarchical coordination strategies that preserve the discrete nature of DER and market interactions.

In recent years, a new strategy has emerged in which the decision order is inverted, so that the discrete operational complexity of active DNs is resolved at the first stage. For example, the authors in \cite{ZHANG2024123073} propose a sequential hierarchical market structure in which virtual power plants first clear a hybrid local market, preserving binary DER scheduling and trading decisions, and the DSO subsequently validates network feasibility through a second-stage flexibility market. Similarly, the work by \cite{10251469} develops a sequential hierarchical coordination framework in which EV charging stations perform P2P scheduling under a fully distributed ADMM negotiation while preserving binary BESS operational decisions. In addition, \cite{10943263} considers a hierarchical setting where multiple DSOs first procure local flexibility and then forward residual bids to the TSO. Although the work addresses grid-safe bid forwarding through corrective and preventive mechanisms, flexibility is limited to continuous active-power adjustments, without PV, BESS, or discrete DER behaviors, despite being one of the few studies examining DSO-first, TSO-second coordination. Along the same line, \cite{garciamunoz2026dso} extends this decision-first paradigm through a bilevel TSO–DSO coordination framework that incorporates PV generation, BESS operation, and P2P energy exchanges within active DNs, showing that DSO-led coordination improves DER utilization relative to conventional TSO-led schemes. Finally, the work in \cite{10720459} analyzes the coordinated provision of voltage regulation and power smoothing in MV–LV distribution grids, where a single MV DSO interacts with a single LV network operator, and the coupling is modeled in OpenDSS with explicit PV and BESS representation. 


The above analysis reveals a growing interest in accurately modeling coordinated decision-making across multiple voltage levels in modern power systems, alongside an evident shift toward more detailed representations of ADNs featuring diverse actors, DER integration, and emerging local market structures. Thus, despite the progress achieved with hierarchical formulations, both two-stage and bilevel, current implementations still face some limitations. In particular, most existing approaches position prosumers or LV-level entities as followers within the hierarchical structure, meaning their operational decisions are determined only after upstream actors have cleared their own schedules. Under this decision ordering, the subsequent application of traditional solution techniques, such as centralized KKT-based reformulations or distributed decomposition methods, typically requires relaxing discrete variables associated with DER scheduling and local market participation. As a result, the representation of flexibility at the distribution level is often simplified, leading to solutions that may lack operational and economic realism. Only in the most recent contributions has a shift emerged in decision ordering, whereby LV entities (or DSO in \cite{garciamunoz2026dso}) with stronger discrete operational characteristics take the lead, and upstream operators subsequently validate or adjust their resources. While this strategy provides an alternative path for preserving discrete behaviors, most studies remain limited to sequential two-stage hierarchies, thereby failing to fully capture strategic interaction and the optimality conditions inherent in bilevel structures. Furthermore, the literature has not yet explored settings in which multiple LV networks act simultaneously as leaders, nor the possibility of enabling energy exchanges across LV grids, mediated by the MV level. Therefore, building on the emerging DSO-first hierarchical trend explored in \cite{garciamunoz2026dso,10943263}, this work addresses the remaining gaps by extending hierarchical coordination from the TSO–DSO interface to the MV–LV interface, through the following contributions:

\begin{itemize}
    \item A bottom-up coordination scheme is introduced to coordinate multiple LV-DNs operating as energy communities with local PV and BESS resources, interconnected through an upstream MV-DN that enables energy exchange among communities. The MV layer ensures operational feasibility while supplying additional energy via conventional generation and shared MV-level PV resources.
    \item A bilevel optimization framework is proposed in which hierarchical coordination among energy communities is structured as a Stackelberg game, where multiple LV networks operate independently in terms of local decision-making, while a market operator aggregates their energy exchange decisions and acts as the leader, and the MV operator acts as the follower, ensuring efficient system feasibility
    \item A framework that preserves the integrality of LV decisions, enabling an exact single-level reformulation of the bilevel problem, in contrast to top-down approaches that require continuous relaxations of discrete variables.
    \item A distributed coordination algorithm based on Lagrangian Dual Decomposition (LDD) with ADMM (LDD-ADMM) is proposed to enable parallel and privacy-enhancing coordination of LV communities by decomposing the coordination structure derived from the bilevel problem. 
\end{itemize}

\begin{table*}[htpb]
\centering
\footnotesize
\renewcommand{\arraystretch}{1.1}
\setlength{\tabcolsep}{3pt}
\begin{tabular}{c c c c c c c c c c}
\hline
\textbf{Ref} & \textbf{Year} & \textbf{Hierarchy} & \textbf{Decision Sequence} & \textbf{OPF} &
\textbf{Resolution} & \textbf{Approach} &
\textbf{PV} & \textbf{BESS} & \textbf{EE} \\
\hline
\cite{YUAN2017600} & 2017 & Two-Stage & TSO$\rightarrow$DSO (1$\rightarrow$n) & $\checkmark$ & Distributed & BD &
$\times$ & $\times$ & $\times$ \\
\cite{9503337} & 2022 & Two-Stage & TSO$\rightarrow$DSO (1$\rightarrow$n) & $\checkmark$ & Distributed & Surrogate LR &
$\times$ & $\times$ & $\times$ \\
\cite{9615006} & 2022 & Two-Stage & TSO$\rightarrow$DSO (1$\rightarrow$n) & $\checkmark$ & Centralized & Solver &
$\checkmark$ & $\times$ & $\times$ \\
\cite{10202840} & 2023 & Two-Stage & TSO$\rightarrow$DSO (1$\rightarrow$n) & $\checkmark$ & Distributed & ADMM &
$\checkmark$ & $\checkmark$ & $\checkmark$ \\
\cite{10124221} & 2024 & Two-Stage & TSO$\rightarrow$DSO (1$\rightarrow$1) & $\checkmark$ & Distributed & CCG-ADMM &
$\checkmark$ & $\checkmark$ & $\times$ \\
\cite{10892311} & 2025 & Bilevel & TSO$\rightarrow$DSO (1$\rightarrow$1) & $\checkmark$ & Distributed & ATC-Q-L &
$\checkmark$ & $\checkmark$ & $\times$ \\
\midrule
\cite{bahramara2015modelling} & 2015 & Bilevel & MO$\rightarrow$LVs (1$\rightarrow$n) & $\times$ & Centralized & KKT &
$\checkmark$ & $\checkmark$ & $\times$ \\
\cite{BAHRAMARA2016169} & 2016 & Bilevel & MO$\rightarrow$LVs (1$\rightarrow$n) & $\times$ & Centralized & KKT &
$\checkmark$ & $\times$ & $\times$ \\
\cite{FERNANDEZ2021123254} & 2021 & Bilevel & MO$\rightarrow$Users (1$\rightarrow$n) & $\times$ & Centralized & KKT &
$\checkmark$ & $\checkmark$ & $\times$ \\
\cite{LECADRE2019317} & 2019 & Bilevel & DSO$\rightarrow$TSO (1$\rightarrow$1) & $\checkmark$ & C\&D & GNE-KKT &
$\checkmark$ & $\times$ & $\times$ \\
\cite{GAZAFROUDI2021117575} & 2021 & Two-Stage & MO$\rightarrow$DSO (1$\rightarrow$1) & $\checkmark$ & Centralized & MILP &
$\checkmark$ & $\times$ & $\checkmark$ \\
\cite{IRIA2022122962} & 2022 & Two-Stage & Agg$\rightarrow$DSO (1$\rightarrow$1) & $\checkmark$ & Distributed & ADMM &
$\checkmark$ & $\checkmark$ & $\times$ \\
\midrule
\cite{9591485} & 2022 & Bilevel & DSO$\rightarrow$EVCS (1$\rightarrow$n) & $\checkmark$ & Distributed & ATC-ADMM &
$\checkmark$ & $\times$ & $\checkmark$ \\
\cite{9645337} & 2022 & Bilevel & MO$\rightarrow$Users (1$\rightarrow$n) & $\times$ & Centralized & KKT &
$\checkmark$ & $\checkmark$ & $\checkmark$ \\
\cite{9721228} & 2023 & Bilevel & MO$\rightarrow$Users (1$\rightarrow$n) & $\checkmark$ & Distributed & Algorithm &
$\checkmark$ & $\checkmark$ & $\checkmark$ \\
\cite{WANG2023109065} & 2023 & Bilevel & DSO$\rightarrow$MGs (1$\rightarrow$n) & $\checkmark$ & Distributed & ADMM-ATC &
$\times$ & $\checkmark$ & $\checkmark$ \\
\cite{JIA2024123259} & 2024 & Bilevel & MO$\rightarrow$Users (1$\rightarrow$n) & $\checkmark$ & Distributed & ADMM &
$\checkmark$ & $\checkmark$ & $\checkmark$ \\
\midrule
\cite{ZHANG2024123073} & 2024 & Two-Stage & VPPs$\rightarrow$DSO (n$\rightarrow$1) & $\checkmark$ & C\&D & BD &
$\checkmark$ & $\checkmark$ & $\checkmark$ \\
\cite{10251469} & 2024 & Two-Stage & EVCS$\rightarrow$DSO (n$\rightarrow$1) & $\checkmark$ & Distributed & ADMM &
$\checkmark$ & $\checkmark$ & $\checkmark$ \\
\cite{10943263} & 2025 & Two-Stage & DSOs$\rightarrow$TSO (n$\rightarrow$1) & $\checkmark$ & Centralized & Algorithm &
$\times$ & $\times$ & $\times$ \\
\cite{10720459} & 2025 & Two-Stage & MV$\rightarrow$LV (1$\rightarrow$1) & $\checkmark$ & Centralized & OpenDSS &
$\checkmark$ & $\checkmark$ & $\times$ \\
\cite{garciamunoz2026dso} & 2026 & Bilevel & DSO$\rightarrow$TSO (n$\rightarrow$1) & $\checkmark$ & Centralized & KKT &
$\checkmark$ & $\checkmark$ & $\times$ \\
This work & 2026 & Bilevel & LVs$\rightarrow$MV (n$\rightarrow$1) &
$\checkmark$ & C\&D & KKT/LDD-ADMM &
$\checkmark$ & $\checkmark$ & $\checkmark$ \\
\hline
\end{tabular}
\caption{ADMM: Alternating Direction Method of Multipliers; Agg: Aggregator; ATC: Analytical target cascading; BD: Benders' decomposition; C\&D: Centralized and Distributed; DSO: Distribution system operator; EVCS: Electric vehicle charging station; EE: Energy Exchange; GNE: Generalized Nash equilibrium; LR: Lagrangian relaxation; LV: Low-voltage; MO: Market operator; MV: Medium voltage; TSO: Transmission system operator;  VPP: Virtual power plant.}
\label{Literature_Review_Table}
\end{table*}
Finally, a structured summary of the most relevant contributions analyzed in this section is provided in Table~\ref{Literature_Review_Table}, allowing a comparison across technical features, hierarchical structures, and coordination schemes. The Table is organized into blocks reflecting how the references are introduced throughout the literature review. The first block groups TSO–DSO coordination schemes; the second collects early distribution-focused models with limited network representation or small-scale settings; the third comprises bilevel approaches with a more active role of DNs, typically under one-to-many decision structures; and the final block highlights recent many-to-one formulations, where multiple downstream agents act first, most commonly within two-stage frameworks. In particular, the \textit{Decision Relation} column identifies the decision-making order (i.e., which actor moves first and which follows) and the level of aggregation involved (e.g., $1{\rightarrow}n$), while the \textit{OPF} column specifies whether network constraints are explicitly modeled (e.g., AC-OPF, linearized OPF, or DC-OPF). Additionally, the \textit{EE} column indicates whether the reviewed works address energy exchanges.


\section{Bottom-Up coordination scheme}\label{sec:Coordination_scheme}
This section introduces the proposed bottom-up coordination scheme, its assumptions and scope, as well as the mathematical models of the interconnected LV and MV-DNs. The main notation used throughout this section and the subsequent formulations is summarized in the nomenclature Table \ref{Nomenclature_Table}.

\subsection{Coordination Scheme}\label{Scheme_Assumptions}

Each LV-DN is modeled as an energy community comprising prosumers equipped with PV or BESS, and traditional consumers. These communities are not fully self-sufficient, and at certain times of day, the upstream MV network must supply or absorb the LVs' energy requirements, creating interdependence between LV and MV decisions, in which the feasibility and cost of LV-level exchanges depend on the MV system's optimal response. In this regard, each community aims to minimize its operating costs, in particular by reducing energy imports from the upstream MV grid, while coordinating local resources and demand. Accordingly, at each operating period every LV-DN submits, at its point of common coupling (PCC), a net power exchange composed of three quantities, whose associated prices are assumed to be known a priori: (i) the cheap import, associated with MV-level PV generation; (ii) the expensive import, supplied by conventional dispatchable generation; and (iii) the exported surplus power offered to the MV-DN when local generation exceeds demand. The resulting net exchange of these three quantities is submitted to a centralized coordinator that aggregates the LV-level decisions and interfaces with the MV-DN, which acts as an operational coordination layer whose primary role is to balance power flows and accommodate the exchange requirements declared by the LV communities while allocating the available generation resources in a technically efficient manner. These resources include utility-scale PV and dispatchable technologies such as cogeneration, whose operational costs are represented by exogenous cost signals from the corresponding asset owners. This structure is based on modern European distribution systems, where the increasing integration of distributed generation has led MV networks to host both renewable and dispatchable resources, supporting local autonomy and operational efficiency. 

\begin{figure}[h]
    \centering
    \includegraphics[scale=0.5]{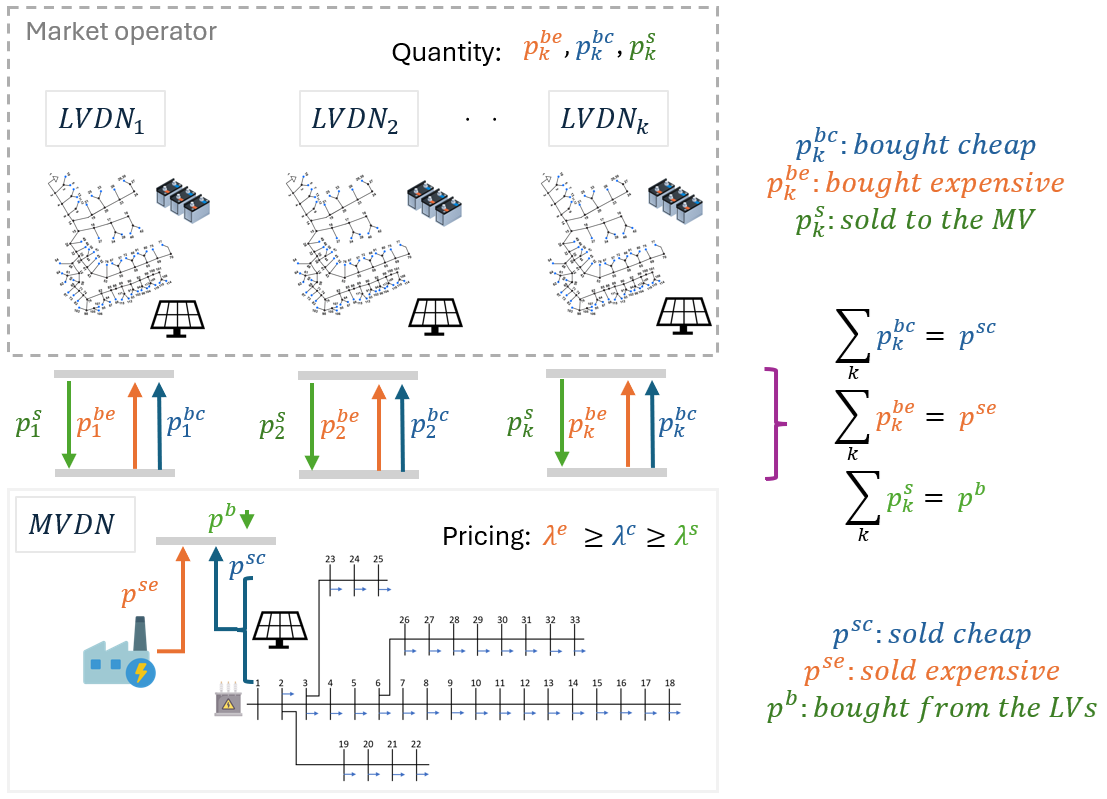}
    \caption{Schematic representation of the bottom-up coordination between LV-DNs and the MV-DN}\label{Esquema_LVs_MV}
\end{figure}

This interdependence between LV exchange decisions and MV operational responses introduces different priorities across voltage levels. While the upper-level coordinator seeks to reduce LV communities' dependence on expensive imports by prioritizing the use of low-cost resources available at the MV level, the MV system generally favors an operational allocation that mitigates congestion, reduces redispatch requirements, and preserves overall network efficiency. Consequently, the MV problem represents an operational response that adapts the allocation of available resources to the exchange requirements imposed by the upper level. As a result, LV decisions may induce operating conditions that differ from the technically efficient allocation preferred by the MV level, thereby giving rise to the strategic interaction captured by the proposed Stackelberg structure. Under this structure, the aggregated LV coordinator acts as the leader by determining the LV communities' exchange requirements, while the MV-DN acts as the follower by optimally responding through the operational allocation of available resources, subject to network constraints. Accordingly, the proposed coordination scheme can be naturally formulated as a bilevel optimization problem, whose hierarchical interaction is schematically illustrated in Figure \ref{Esquema_LVs_MV}.


The above coordination scheme follows a bottom-up decision structure, in which entities at lower layers of the traditional power-system hierarchy make operational decisions before the upstream MV level responds. This inversion of the conventional top-down coordination paradigm is motivated by the evolving decision-making structure of modern DNs and their increasing interaction with upstream systems. As discussed in the literature, the growing penetration of DERs, such as BESS, EV charging, peer-to-peer trading, and flexibility activation mechanisms, introduces an increasing number of discrete and non-convex decision variables at the distribution level. Under conventional top-down coordination schemes, these decisions are typically embedded in the lower-level problem, yielding a non-convex follower that violates the assumptions underlying standard bilevel solution approaches, including strong duality and KKT-based reformulations. In contrast, the proposed bottom-up structure assigns the discrete operational decisions to the upper level while preserving a continuous and convex representation of the MV-DN at the lower level. This structural separation ensures that the follower problem satisfies convexity, primal feasibility, and strong duality conditions, thereby enabling an exact reformulation of the proposed bilevel optimization problem.


Unlike conventional bilevel formulations, where prices are typically determined endogenously through market-clearing mechanisms or dual variables, the exchange prices considered in this work are treated as exogenous parameters. This assumption is based on the operational nature of the proposed coordination framework and on general current DN practices, where the economic signals faced by LV communities are commonly defined by externally established tariff structures, bilateral agreements, or regulated remuneration schemes rather than by endogenous market interactions. In particular, the selling price of exported surplus energy reflects a predefined feed-in remuneration mechanism, while the import prices of cheap and expensive energy correspond to operational price signals linked to the expected availability and operational costs of the resources connected at the MV level. These prices are fixed before the LV exchange decisions are made and are not obtained from the lower-level MV optimization problem. Nevertheless, this assumption does not preclude the economic interpretation of the dual variables associated with the follower problem, whose marginal values may still provide information regarding congestion conditions and the marginal operational cost of allocating cheap and expensive resources within the MV network.

Within the scope of the proposed coordination framework, the MV-DN is modeled as a self-contained system, with no energy exchange with the upstream TN. In other words, no power is allowed to enter or leave the MV level beyond the resources explicitly represented in the model. System feasibility is ensured exclusively by dispatchable conventional generation and utility-scale PV units integrated at the MV level.  Accordingly, the LV exchange schedules are assumed to be admissible only when they can be accommodated by the MV-DN without requiring upstream imports or exports. Likewise, surplus energy exported by the LV communities can always be absorbed within the MV-DN. This is guaranteed by the presence of conventional MV demand nodes that can accommodate redistributed surplus generation, thereby preventing the need for upstream export. This modeling choice avoids introducing additional layers of complexity associated with upstream market interactions or transmission-level coupling, which are not essential for the objectives pursued in this work. By restricting the system to the explicitly modeled MV and LV levels, the analysis can concentrate on the structural properties of the proposed bottom-up coordination scheme and on the performance of the exact and distributed resolution strategies. The emphasis is therefore placed on validating the bilevel formulation and its algorithmic implementation, while extensions that incorporate broader system boundaries and upstream integration are left for future research.

\subsection{Low voltage DN model}\label{LV_DN_Model}

Consider a set of LV-DNs indexed by $k \in \mathcal{K}$, each connected to a subset of MV buses $\Theta^{k} \subseteq \Theta$. Each LV-DN $k$ is represented as a radial directed graph $\mathcal{G}_k = (\Omega_k, \mathcal{L}_k)$, where $\Omega_k$ denotes the set of buses and $\mathcal{L}_k$ the set of distribution lines, and is operated over a finite scheduling horizon $t \in \mathcal{T}$ with a fixed time resolution $\Delta t$. A single bus $i \in \Omega_k^{B} \subseteq \Omega_k$ is designated as the PCC, through which the LV-DN interfaces electrically with the upstream MV network. The internal operation of each LV-DN is formulated as a mixed-integer second-order cone programming (MISOCP) problem based on the DistFlow model \cite{19266}. The formulation explicitly captures active and reactive power flows, nodal voltage magnitudes, and line current limits, while preserving convexity through a second-order cone relaxation of the AC power flow equations \cite{6756976}. The set of buses $\Omega_k$ is partitioned into different subsets. The subset $\Omega_k^{A}\subseteq \Omega_k$ corresponds to user buses that may be equipped with DERs, including PV generation units and/or BESS, while the remaining buses represent non-user internal network nodes without controllable resources. Accordingly, the objective function minimizes the total operational cost over the scheduling horizon, which is formulated as follows:

\begin{align} \label{OF_LV_DN_k}
    \min \sum_{t \in \mathcal{T}} \sum_{i \in \Omega^B_k} & \Bigg[
         \left( \lambda^{c}_{t}  Mp^{bc}_{i,t}
        + \lambda^{e}_{t}  Mp^{be}_{i,t}
        - \lambda^{s}_{t}  Lp^{s}_{i,t} \right) \Bigg] + 
         \sum_{t \in \mathcal{T}} \sum_{i \in \Omega^A_k} \Bigg[ C^{pv}  pv^{\mathrm{L}}_{i,t} + C^{bt}  ds^{\mathrm{L}}_{i,t}
    \Bigg]. 
\end{align}

The objective function in \eqref{OF_LV_DN_k} is composed of two main components. The first term represents the net cost of energy exchanges at the PCC. Thus, when the LV-DN cannot fully meet its internal demand using local resources, it imports energy from the upstream MV network, considering cheap energy ($Mp^{bc}_{i,t}$), associated with utility-scale PV generation, and expensive energy ($Mp^{be}_{i,t}$) from conventional dispatchable generation. Conversely, when local generation exceeds demand, the LV-DN exports surplus energy ($Lp^{s}_{i,t}$) to the MV network, which is reflected in the objective function as a revenue term. These boundary exchanges are valued using time-dependent prices ($\lambda^{(\cdot)}_{t}$), assumed to be known a priori and treated as exogenous parameters within the coordination framework. The second component of the objective function represents the internal operational costs associated with the use of DERs at the user level. This includes the cost of PV generation ($C^{pv}$), which may reflect feed-in tariff schemes or technology-specific operational costs, and the cost of BESS discharging ($C^{bt}$). The latter is introduced to capture the effects of battery degradation and represents the implicit depreciation cost incurred when using stored energy. Accordingly, the minimization of \eqref{OF_LV_DN_k} is subject to the network and DER operational constraints introduced below.

\begin{subequations}\label{LV_operation_Eqs}
\begin{align}
    &\displaystyle\sum_{(i,j) \in {\mathcal{L}_k}} p^{\mathrm{L}}_{i,j,t} - \displaystyle\sum_{(j,i) \in {\mathcal{L}_k}} (p^{\mathrm{L}}_{j,i,t} - R^{\mathrm{L}}_{j,i}\ell^{\mathrm{L}}_{i,j,t}) =  
    \left \{
      \begin{aligned}
            & \Delta p^{\mathrm{L}}_{i,t} && \text{if } \forall i \in \Omega^A_k\\
            & \Delta b^{\mathrm{L}}_{i,t} && \text{if } \forall i \in \Omega^B_k\\
        & 0 && Otherwise
      \end{aligned} \right. && \forall i \in \Omega_k, \forall t \in \mathcal{T}\label{LV_OP_Eq_PG} \\
    &\displaystyle\sum_{(i,j) \in {\mathcal{L}_k}} q^{\mathrm{L}}_{i,j,t} - \displaystyle\sum_{(j,i) \in {\mathcal{L}_k}} (q^{\mathrm{L}}_{j,i,t} - X^{\mathrm{L}}_{j,i}\ell^{\mathrm{L}}_{i,j,t}) =  qg^{\mathrm{L}}_{i,t} - QL^{\mathrm{L}}_{i,t} && \forall i \in \Omega_k, \forall t \in \mathcal{T} \label{LV_OP_Eq_QG}  \\
    &\Delta p^{\mathrm{L}}_{i,t} = pv^{\mathrm{L}}_{i,t} - PL^{\mathrm{L}}_{i,t} + ds^{\mathrm{L}}_{i,t} - ch^{\mathrm{L}}_{i,t} && \forall i \in \Omega^A_k, \forall t \in \mathcal{T} \label{LV_OP_Eq_DP} \\
    &\Delta b^{\mathrm{L}}_{i,t} = Mp^{bc}_{i,t} + Mp^{be}_{i,t} - Lp^{s}_{i,t}  && \forall i \in \Omega^B_k, \forall t \in \mathcal{T} \label{LV_OP_Eq_DB} \\
    &v^{\mathrm{L}}_{j,t} = v^{\mathrm{L}}_{i,t} -2 (R^{\mathrm{L}}_{i,j} p^{\mathrm{L}}_{i,j,t} + X^{\mathrm{L}}_{i,j}q^{\mathrm{L}}_{i,j,t}) + ({R^{\mathrm{L}}}^2_{i,j} + {X^{\mathrm{L}}}^2_{i,j})\ell^{\mathrm{L}}_{i,j,t} && \forall (i,j) \in \mathcal{L}_k, \forall t \in \mathcal{T} \label{LV_OP_Eq_V} \\
    &(p^{\mathrm{L}}_{i,j,t})^2 + (q^{\mathrm{L}}_{i,j,t})^2 \leq \ell^{\mathrm{L}}_{i,j,t} v^{\mathrm{L}}_{i,t} && \forall (i,j) \in \mathcal{L}_k, \forall t \in \mathcal{T} \label{LV_OP_Eq_S}\\
    &Q^{\mathrm{L}}_{min} \leq qg^{\mathrm{L}}_{i,t} \leq Q^{\mathrm{L}}_{max} && \forall i \in \Omega_k, \forall t \in \mathcal{T} \label{LV_OP_Eq_Q}\\
    &{V^{\mathrm{L}}}^{min}_{i} \leq v^{\mathrm{L}}_{i,t} \leq {V^{\mathrm{L}}}^{max}_{i}  && \forall i \in \Omega_k, \forall t \in \mathcal{T}  \label{LV_OP_Eq_VL}\\
    &\ell^{\mathrm{L}}_{i,j,t} \leq {I^{\mathrm{L}}}^{max}_{i,j} && \forall i \in \Omega_k, \forall t \in \mathcal{T} \label{LV_OP_Eq_ell}
\end{align}
\end{subequations}

Eqs.~\eqref{LV_operation_Eqs} describe the network operation of each LV-DN following a SOCP relaxation of the DistFlow model. Specifically, Eqs.~\eqref{LV_OP_Eq_PG} and \eqref{LV_OP_Eq_QG} represent the nodal active and reactive power balance constraints, respectively. For user buses, the active power imbalance in Eq.~\eqref{LV_OP_Eq_PG} is defined by the net injection $\Delta p^{L}_{i,t}$, which is detailed in Eq.~\eqref{LV_OP_Eq_DP} as the balance between local PV generation ($pv^{\mathrm{L}}_{i,t}$), electric demand ($PL^{\mathrm{L}}_{i,t}$), and BESS charging ($ch^{\mathrm{L}}_{i,t}$) and discharging ($ds^{\mathrm{L}}_{i,t}$). At the PCC, the net boundary exchange $\Delta b^{L}_{i,t}$ is defined in Eq.~\eqref{LV_OP_Eq_DB} as the difference between imported cheap and expensive energy and exported surplus power. Eq.~\eqref{LV_OP_Eq_V} models the voltage drop along distribution lines, while Eq.~\eqref{LV_OP_Eq_S} corresponds to the second-order cone constraint linking active and reactive power flows with squared current magnitudes. Finally, Eqs.~\eqref{LV_OP_Eq_Q}–\eqref{LV_OP_Eq_ell} enforce operational limits on reactive power exchange, nodal voltage magnitudes, and line current capacities, respectively. Note that the superscript $(\cdot)^{L}$ is used throughout the formulation to
differentiate variables, constraints, and parameters associated with the LV-DN from their counterparts defined at the MV level. The operation of PV units and BESS is subsequently modeled using the following set of constraints.

\begin{subequations}\label{LV_DER_Eqs}
    \begin{align}
        &pv^{\mathrm{L}}_{i,t} \leq {\Gamma^{\mathrm{L}}}^{pv}_{i} {PV^{\mathrm{L}}}_{t} && \forall i \in \Omega^A_k, \forall t \in \mathcal{T} \label{LV_OP_Eq_PV}\\
        &soc^{\mathrm{L}}_{i,t}= soc^{\mathrm{L}}_{i,t-1}+\left(\varphi^{ch} ch^{\mathrm{L}}_{i,t}-\frac{1}       {\varphi^{ds}}ds^{\mathrm{L}}_{i,t}\right)\Delta t && \forall i \in \Omega^A_k, \forall t \in \mathcal{T} \label{LV_OP_Eq_SOC} \\
        &{\Gamma^{\mathrm{L}}}^{bt}_{i} {SOC^{\mathrm{L}}}^{min}\leq soc^{\mathrm{L}}_{i,t} \leq {\Gamma^{\mathrm{L}}}^{bt}_{i} {SOC^{\mathrm{L}}}^{max} && \forall i \in \Omega^A_k, \forall t \in \mathcal{T} \label{LV_OP_Eq_SOCL} \\
        &ch^{\mathrm{L}}_{i,t}  \leq {PB^{\mathrm{L}}}_{i} (1 - w^{\mathrm{L}}_{i,t}) - {PB^{\mathrm{L}}}_{i} (1- \nu_i) && \forall i \in \Omega^A_k, \forall t \in \mathcal{T} \label{LV_OP_Eq_CH}\\
        &ds^{\mathrm{L}}_{i,t} \leq {PB^{\mathrm{L}}}_{i} w^{\mathrm{L}}_{i,t} && \forall i \in \Omega^A_k, \forall t \in \mathcal{T} \label{LV_OP_Eq_DS}\\
        &w^{\mathrm{L}}_{i,t} \leq \nu_i && \forall i \in \Omega^A_k, \forall t \in \mathcal{T} \label{LV_OP_Eq_W}
    \end{align}
\end{subequations}

Eqs.~\eqref{LV_DER_Eqs} model the operation of PV generation units and BESS at user buses. Specifically, Eq.~\eqref{LV_OP_Eq_PV} limits the PV power injection as a function of the installed capacity, where $\Gamma^{L}_{pv}$ denotes the PV capacity and $PV^{L}_{t}$ represents the normalized solar irradiance profile. Eqs.~\eqref{LV_OP_Eq_SOC}–\eqref{LV_OP_Eq_W} describe the BESS operation. Eq.~\eqref{LV_OP_Eq_SOC} models the state-of-charge dynamics, while Eq.~\eqref{LV_OP_Eq_SOCL} enforces minimum and maximum SOC limits, defined as percentages of the installed BESS capacity $\Gamma^{L}_{bt}$. Charging and discharging decisions are governed by Eqs.~\eqref{LV_OP_Eq_CH} and \eqref{LV_OP_Eq_DS}, which depend on the binary variable $w^{L}_{i,t}$ to prevent simultaneous charging and discharging, and on the parameter $\nu_i$ indicating the presence of a BESS at the corresponding node in Eq.~\eqref{LV_OP_Eq_W}. Finally, $PB^{L}_{i}$ represents the maximum charging and discharging power of the BESS.

\subsection{Medium voltage DN model}\label{MV_DN_Model}
Similar to the LV-DNs, the MV-DN is modeled using a DistFlow-based formulation and is represented as a directed graph $\mathcal{G} = (\Theta, \mathcal{L})$, where $\Theta$ denotes the set of MV buses and $\mathcal{L}$ the set of distribution lines. The MV-DN is operated over the same scheduling horizon $t \in \mathcal{T}$ and time resolution $\Delta t$ to ensure temporal consistency with the LV-level decisions and system feasibility. As discussed in Section~\ref{Scheme_Assumptions}, the MV-DN acts as a coordinating layer that balances power flows across the network and supplies the aggregated demands declared by the connected LV-DNs at minimum operational cost. In this regard, the set of MV buses $\Theta$ is partitioned into different subsets. The subset $\Theta^{k} \subseteq \Theta $ corresponds to the buses at which LV-DNs are connected to the MV network. The subset $\Theta^{A} \subseteq \Theta $ represents internal MV buses that may host conventional demand, modeled through fixed demand profiles, or be connected to utility-scale resources such as PV generation units. Accordingly, the MV-DN optimization problem determines the dispatch of conventional generation, utility-scale PV to satisfy the aggregated power requirements of the LV-DNs and the internal MV demand, while minimizing the system's total operating cost. This objective is formulated as a convex quadratic economic dispatch problem, as defined below.
\begin{align}\label{MV_OF}
    \min \quad & \sum_{t \in \mathcal{T}} \sum_{i \in \Theta^A} \left[
        C^{a}_{i} (pg_{i,t})^2
        + C^{b}_{i} (pg_{i,t})
        + C^{c}_{i}
        + C^{pv} pv^M_{i,t}
    \right] 
\end{align}

The objective function in Eq.~\eqref{MV_OF} represents a convex quadratic economic dispatch problem at the MV level. Thus, the first three terms correspond to the generation cost of conventional thermal units, where $pg^{M}_{i,t}$ denotes the dispatched active power, and $C^{a}_{i}$, $C^{b}_{i}$, and $C^{c}_{i}$ represent the quadratic, linear, and fixed
cost coefficients, respectively. In addition, the objective function in Eq.~\eqref{MV_OF} includes operational cost terms associated with the utilization of utility-scale PV generation modeled in an analogous manner to the LV-DN formulation. Note that the superscript $(\cdot)^{M}$ is used to distinguish variables and parameters associated with the MV-DN from their counterparts defined at the LV level. Therefore, the objective function is minimized subject to the MV-DN operational and network constraints detailed below.

\begin{subequations} \label{MV_Op_Eqs}
    \begin{align}
        &\displaystyle\sum_{(i,j) \in {\mathcal{L}}} p^{\mathrm{M}}_{i,j,t} - \displaystyle\sum_{(j,i) \in {\mathcal{L}}} p^{\mathrm{M}}_{j,i,t} =  
        \left \{
          \begin{aligned}
                & \Delta p^{\mathrm{M}}_{i,t} && \text{if } \forall i \in \Theta^A\\
                & \Delta b^{\mathrm{M}}_{i,t} && \text{if } \forall i \in \Theta^k\\
            & 0 && Otherwise
          \end{aligned} \right. && \forall i \in \Theta, \forall t \in \mathcal{T} \label{MV_Op_Eqs_PG} \\
        &\displaystyle\sum_{(i,j) \in {\mathcal{L}}} q^{\mathrm{M}}_{i,j,t} - \displaystyle\sum_{(j,i) \in {\mathcal{L}}} q^{\mathrm{M}}_{j,i,t} =  qg^{\mathrm{M}}_{i,t} - QL^{\mathrm{M}}_{i,t} && \forall i \in \Theta, \forall t \in \mathcal{T} \label{MV_Op_Eqs_QG}  \\
        &\Delta p^{\mathrm{M}}_{i,t} = pv^{\mathrm{M}}_{i,t} + pg^{\mathrm{M}}_{i,t} - PL^{\mathrm{M}}_{i,t}  && \forall i \in \Theta^A, \forall t \in \mathcal{T} \label{MV_Op_Eqs_DP}  \\
        &\Delta b^{\mathrm{M}}_{i,t} = - \Delta b^{\mathrm{L}}_{i,t}  && \forall k \in \mathcal{K}, \forall i \in \Theta^k, \forall t \in \mathcal{T} \label{MV_Op_Eqs_DB} \\ 
        &v^{\mathrm{M}}_{j,t} = v^{\mathrm{M}}_{i,t} -2 (R^{\mathrm{M}}_{i,j} p^{\mathrm{M}}_{i,j,t} + X^{\mathrm{M}}_{i,j}q^{\mathrm{M}}_{i,j,t})  && \forall (i,j) \in \mathcal{L}, \forall t \in \mathcal{T} \label{MV_Op_Eqs_V} \\
        &Q^{\mathrm{M}}_{min} \leq qg^{\mathrm{M}}_{i,t} \leq Q^{\mathrm{M}}_{max} && \forall i \in \Theta, \forall t \in \mathcal{T} \label{MV_Op_Eqs_Q} \\
        &{V^{\mathrm{M}}}^{min}_{i} \leq v^{\mathrm{M}}_{i,t} \leq {V^{\mathrm{M}}}^{max}_{i}  && \forall i \in \Theta, \forall t \in \mathcal{T}   \label{MV_Op_Eqs_VL} 
    \end{align}
\end{subequations}
Eqs.~\eqref{MV_Op_Eqs} describe the MV-DN network operation and largely mirror the formulation adopted for the LV-DNs, with a few key differences. In particular, the MV-DN is modeled using a linearized DistFlow representation \cite{9268988}, where network losses are neglected. As a result, the quadratic loss-related terms present in the LV-DN formulation, namely the $R\,\ell$ and $X\,\ell$ components, are omitted, and no SOCP relaxation is required at this level \cite{7484653}. Active and reactive power balance constraints are enforced through Eqs.~\eqref{MV_Op_Eqs_PG} and \eqref{MV_Op_Eqs_QG}, while Eqs.~\eqref{MV_Op_Eqs_DP} and \eqref{MV_Op_Eqs_DB} define the net active power injections at internal MV buses and at the buses where LV-DNs are connected, respectively. Note that, at the PCC buses, the net boundary exchange between the MV-DN and each LV community appears with opposite sign. This sign convention reflects the direction of power flows in the nodal balance: power injected into the MV network is considered positive, whereas power withdrawn from it is negative. Consequently, when an LV community imports power, the MV-DN supplies it (negative injection at the MV node), and conversely, when an LV community exports surplus energy, it is represented as a positive injection into the MV network. Likewise, the operational constraints of utility-scale PV generation units at the MV level follow a structure analogous to that adopted for the LV-DNs. Note that, although network losses are neglected in the linearized formulation, the voltage drop equation~\eqref{MV_Op_Eqs_V} retains the resistance and reactance parameters of each line, meaning that the physical topology of the MV network still influences the nodal voltage profiles and, consequently, the feasible allocation of generation resources. This topological dependence introduces a spatial dimension to the follower's operational response: the cost-minimizing dispatch of cheap and expensive resources is determined not only by the energy balance but also by the electrical locations of LV communities and generation units within the network. As a result, when the follower optimizes its dispatch to accommodate the leader's exchange requirements, the network topology may prevent it from reaching its independently optimal allocation, forcing the MV system to operate at a higher cost than it would under an unconstrained economic dispatch.
\begin{subequations}\label{MV_PG_Eq}
    \begin{align}
        &pv^{\mathrm{M}}_{i,t} \leq {\Gamma^{\mathrm{M}}}^{pv}_{i} {PV_{t}^{\mathrm{M}}} && \forall i \in \Theta^A, \forall t \in \mathcal{T} \label{MV_PG_Eq_PV}\\
        &pg^{\mathrm{M}}_{i,t} \leq {PG_{i,t} ^{\mathrm{M}}} && \forall i \in \Theta^A, \forall t \in \mathcal{T} \label{MV_PG_Eq_Pg}\\
        &\displaystyle\sum_{i \in \Theta^k} Mp^{bc}_{i,t} \leq \displaystyle\sum_{i \in \Theta^A} pv^{\mathrm{M}}_{i,t} && \forall t \in \mathcal{T} \label{MV_PG_Eq_PVM}\\
        &\displaystyle\sum_{i \in \Theta^k} Mp^{be}_{i,t} \leq \displaystyle\sum_{i \in \Theta^A} pg^{\mathrm{M}}_{i,t} &&   \forall t \in \mathcal{T}\label{MV_PG_Eq_PGM}
    \end{align}
\end{subequations}
Eqs.~\eqref{MV_PG_Eq} complete the MV-DN formulation by defining the operation of utility-scale generation resources. Eq.~\eqref{MV_PG_Eq_PV} limits the active power injection from utility-scale PV generation, while Eq.~\eqref{MV_PG_Eq_Pg} corresponds to the dispatchable conventional generation, where $PG^{M}_{i}$ defines the maximum available generation capacity. Eqs.~\eqref{MV_PG_Eq_PVM} and \eqref{MV_PG_Eq_PGM} enforce global balance constraints ensuring that the total amount of cheap and expensive energy supplied to the connected LV-DNs does not exceed the available utility-scale PV generation and conventional generation capacity, respectively. 

\subsection{Bilevel Coordination Problem}
To facilitate the explanation of the bilevel coordination problem formulation, the optimization models for the LV and MV networks are expressed in a semi-compact form that preserves the hierarchical coupling structure and explicitly distinguishes boundary variables. This representation also provides the basis for the distributed decomposition developed in the following section.
\begin{equation}
\label{LV-DNk}
\text{(LV-DN$_k$)} \quad
\left\{
\begin{aligned}
\min_{x_k,\, b_k} \quad
& F_k(x_k)
+ \sum_{t \in \mathcal{T}} \sum_{i \in \Omega_k^B}
\left(
\lambda_t^c Mp^{bc}_{i,t}
+ \lambda_t^e Mp^{be}_{i,t}
- \lambda_t^s Lp^{s}_{i,t}
\right) \\[0.5em]
\text{s.t.} \quad
& G_k(x_k) = 0, \\ 
& H_k(x_k) \le 0, \\
& \Delta b^{L}_{i,t}
= Mp^{bc}_{i,t}
+ Mp^{be}_{i,t}
- Lp^{s}_{i,t},
\quad \forall i \in \Omega_k^B,\ \forall t \in \mathcal{T}
\end{aligned}
\right.
\end{equation}

\noindent
where $x_k$ gathers the internal operational variables of LV-DN$_k$ (voltages, flows, PV, and BESS), and  $\Delta b^{L}_k = \{Mp^{bc}_{k,t}, Mp^{be}_{k,t}, Lp^{s}_{k,t}\}$ represents the boundary exchange variables.  The function $F_k(\cdot)$ denotes the local operational costs associated with generation and storage. 
\begin{equation}
\label{MV-DN}
\text{(MV-DN)} \quad
\left\{
\begin{aligned}
\min_{y,\, b} \quad
& f(y) \\[0.5em]
\text{s.t.} \quad
& g(y) = 0, \\
& h(y) \le 0, \\
& \Delta b^{M}_{i,t}
= -\Delta b^{L}_{i,t} 
\quad &&\forall i \in \Theta^k, \forall t \in \mathcal{T}, 
\end{aligned}
\right.
\end{equation}

\noindent
Here, $y$ represents the internal MV network variables (flows, voltages, generation, and storage),  and $\Delta b_{i,t}^{M} = -\Delta b^{L}_{i,t}$ are the coupling variables between the MV and LV networks. The function $f(y)$ corresponds to the total operational cost of the MV-DN, including centralized PV and traditional generation resources. 

Based on the compact form of the LV-DNs and MV-DN, the bilevel optimization model can be expressed as follows.
\begin{equation}\label{BL_problem}
\text{(BL)} \quad
\left\{
\begin{aligned}
\min_{\{x_k,\Delta b_k^L\}_{k \in \mathcal{K}}} \quad
& \sum_{k \in \mathcal{K}} F_k(x_k)
+ \sum_{k \in \mathcal{K}}\sum_{i \in \Omega_k^B}\sum_{t\in \mathcal{T}}
\left(
\lambda_t^c Mp^{bc}_{i,t}
+ \lambda_t^e Mp^{be}_{i,t}
- \lambda_t^s Lp^{s}_{i,t}
\right) \\[0.6em]
\text{s.t.} \quad
& G_k(x_k) = 0,  \\
& H_k(x_k) \le 0,  \\
& \Delta b^L_{i,t}
= Mp^{bc}_{i,t} + Mp^{be}_{i,t} - Lp^{s}_{i,t},
\quad \forall i \in \Omega_k^B,\ \forall t \in\mathcal{T}\\[0.6em]
& (y, Mp^{(\cdot)} ) \in \arg\min_{y, Mp^{(\cdot)}}
\Big\{
f(y) : g(y) = 0,\ h(y) \le 0, \\
& \qquad\qquad \Delta b^{M}_{i,t} = -\Delta b^{L}_{i,t} , \\
& \qquad\qquad
\forall k \in \mathcal{K}, \forall i \in \Theta^k,\forall t \in \mathcal{T}
\Big\}
\end{aligned}
\right.
\end{equation}

In this formulation \eqref{BL_problem}, the upper level represents the coordinated operation of multiple LV energy communities indexed by $k \in \mathcal K$. Thus, each community $k$ minimizes its operational cost subject to its local operational constraints, $G_k(x_k)=0$ and $H_k(x_k)\le 0$, and to the boundary balance equations defining the net exchange $\Delta b^L_{i,t}$. The decisions of all LV communities are collected by a centralized coordinator that interfaces with the lower-level entities represented by the MV-DN. In this hierarchical interaction, the LV-level coordinator serves as the leader by declaring the boundary exchanges, while the MV-DN serves as the follower by responding optimally through the cost-efficient allocation of cheap and expensive generation resources to ensure system feasibility. Given the boundary exchange quantities, the MV-DN solves an operational optimization problem to minimize the MV-DN cost $f(y)$ subject to the network and operational constraints $g(y)=0$ and $h(y)\le 0$, with both levels coupled through consistency constraints equating the MV-level exchange variables $-\Delta b^{L}_{i,t}$ to their LV-level counterparts. This leader--follower interaction is interpreted as a Stackelberg structure, in which the leader anticipates the MV-DN's optimal operational response and determines its exchange decisions accordingly. In this regard, the proposed formulation adopts the standard optimistic interpretation of bilevel programming, in which, if multiple optimal responses exist for the MV-DN problem, the follower is assumed to select the solution most favorable to the leader.

\section{Distributed coordination algorithm via Lagrangian Relaxation and ADMM}\label{sec:Resolution_strategy}
The bottom-up coordination scheme defined as a BL problem in \eqref{BL_problem} can be reformulated as a single-level (SL) problem by replacing the lower-level problem with its commonly used KKT optimality conditions (see details in \ref{Appendix_A}). This reformulation is exact provided that the lower-level problem is convex, admits an optimal solution, and satisfies strong duality and the corresponding regularity conditions \cite{dempe2002foundations}. Under these conditions, the KKT system is necessary and sufficient for optimality; therefore, the resulting SL formulation is equivalent to the original BL problem and serves as a theoretical benchmark. The complementarity conditions arising from the KKT reformulation can be handled through standard mixed-integer implementations, such as Big-M linearization or SOS1 \cite{aussel2025tutorial} constraints. The latter offers a convenient alternative by enforcing complementarity without requiring user-defined Big-M parameters, thereby mitigating potential numerical issues associated with poorly chosen bounds, which is the approach adopted in this work. However, this SL problem may raise practical concerns because it assumes the existence of a centralized coordinator at the upper level with full access to the operational decisions of all LV-DNs, including detailed information on local PV generation and BESS operation. While this assumption may be acceptable from a theoretical standpoint, it is often undesirable in practice due to concerns regarding data confidentiality and privacy concerns. Preserving the autonomy of individual LV communities and limiting the disclosure of sensitive operational information have been recognized as key requirements in hierarchical coordination schemes for ADN \cite{RANJBAR2024122823}.

Among the available distributed optimization techniques, ADMM has emerged as an effective tool for coordinating decomposable problems with limited information exchange \cite{neal2011distributed}. In this regard, the authors in \cite{nie2025augmented} have proposed the Augmented Lagrangian techniques as a mechanism for reformulating bilevel problems and handling lower-level optimality conditions, and similarly in \cite{ZHAO2024122731}, while ADMM-based distributed schemes have been successfully applied to large-scale coordination problems by exploiting separable structures and reformulating coupling constraints into tractable consensus forms \cite{ZHAO2024122731}. However, the bilevel structure considered in this work cannot be directly addressed through a standard ADMM implementation, as the interaction between LV- and MV-level variables involves coupled decisions that are not readily separable. In particular, a direct relaxation of the aggregate boundary balance constraint, $\Delta b^{M} = -\Delta b^{L}$, is insufficient to ensure stable coordination. This limitation arises because the net exchange variable $\Delta b$ aggregates multiple underlying transactions (cheap imports, expensive imports, and surplus exports), thereby introducing degrees of freedom that are not uniquely identified when only the net balance is coordinated. As a result, different combinations of these components may yield the same net exchange, leading to poor convergence behavior. Therefore, an additional reformulation is required to restructure the coupling and enable a distributed solution scheme.

To address this limitation, the aggregate boundary exchange is decomposed into its underlying components, and auxiliary copies of these variables are introduced to represent the boundary decisions at the LV and MV sides  ($ Lp^s_{i,t} = Mp^s_{k,t}, \enspace Mp^{be}_{k,t} = Lp^{be}_{i,t}, \enspace Mp^{bc}_{k,t} = Lp^{bc}_{i,t}$),  where $Lp^{bc}_{i,t}, Lp^{be}_{i,t}, Mp^s_{k,t}$ denote the auxiliary variables.  The LV and MV problems are then rewritten in terms of their respective local copies, while consistency between these variables is enforced through coupling constraints defined outside the individual subproblems \cite{tosserams2007augmented} (see Appendix \ref{Appendix_B}).  These constraints do not belong to either the leader or the follower individually, but instead define the interaction between the two levels in the reformulated problem.  This type of reformulation, based on variable duplication and consensus constraints, is adopted by \cite{wei2026low} to enable ADMM-based solution schemes in problems with coupled decision structures.  However, enforcing all three consistency constraints is not necessary.  By preserving the original net balance relation, consistency can be imposed on only two components, with the third implicitly determined.


Based on this reformulated coupling structure (DC-BL \eqref{BL_reformulated_full}), a hybrid coordination scheme can be adopted. The introduction of duplicate boundary variables, together with consistency constraints, prevents a direct separable solution to the coordination problem, as these constraints couple the LV and MV decisions. To enable decomposition, these coupling constraints are relaxed through a Lagrangian framework. In particular, the coordination problem is expressed through a partially augmented Lagrangian, in which the consistency constraints associated with the duplicated boundary variables are dualized, while only a subset of them is reinforced via quadratic penalty terms. Specifically, the consistency constraint associated with the cheap-energy component is incorporated into the augmented Lagrangian, introducing quadratic penalty terms that couple the LV and MV copies of this variable. In contrast, the surplus exchange constraint is handled through standard Lagrangian relaxation, entering the Lagrangian only through its corresponding dual multiplier without additional penalization. The remaining component, associated with expensive energy imports, is implicitly recovered from the net boundary balance, and therefore does not require explicit dualization.

The resulting Lagrangian function is then minimized with respect to the primal variables. Due to the separable structure induced by the variable duplication and the partial augmentation of the Lagrangian, this minimization can be decomposed into independent blocks corresponding to each LV-DN and the MV-DN. This decomposition yields a sequence of coordinated subproblems, where each block inherits a specific objective function composed of its original cost terms, augmented by the corresponding dual variables, and, in the case of the cheap-energy component, quadratic penalty terms that enforce consensus. In particular, the LV subproblems include the augmented penalty for deviations of their local cheap-energy variables from the shared consensus value, while the MV subproblem incorporates the corresponding dual and penalty terms reflecting its role in allocating this resource, as well as the dualized surplus exchange. In this way, the proposed LDD-ADMM algorithm minimizes the hybrid Lagrangian via a block-wise decomposition that preserves the bilevel reformulation's behavior, reduces the dimensionality of the consensus space, and yields a sequence of independent LV- and MV-level subproblems that can be solved in parallel.
\begin{equation}\label{LV_ADMM_problem}
\text{(LV-DN$_k$--ADMM)} \quad
\left\{
\begin{aligned}
\min_{x_k,\, \Delta b_k^L} \quad
& F_k(x_k)
+ \sum_{t \in \mathcal{T}} \sum_{i \in \Omega_k^B}
\left(
\lambda_t^{c} Lp^{bc}_{i,t}
+ \lambda_t^e Lp^{be}_{i,t}
- \lambda_t^s Lp^{s}_{i,t}
\right) \\
& \quad
+ \sum_{t \in \mathcal{T}} \sum_{i \in \Omega_k^B}
\frac{\rho}{2}
\left\|
Lp^{bc}_{i,t}
- Mp^{bc,(r)}_{k,t}
+ u^{bc,(r)}_{k,t}
\right\|^2 \\[0.5em]
\text{s.t.} \quad
& G_k(x_k) = 0, \\
& H_k(x_k) \le 0, \\
& \Delta b^{L}_{i,t}
= Lp^{bc}_{i,t}
+ Lp^{be}_{i,t}
- Lp^{s}_{i,t},
\quad
\forall i \in \Omega_k^B,\ \forall t \in\mathcal{T}.
\end{aligned}
\right.
\end{equation}

\begin{equation}\label{MV_ADMM_problem}
\text{(MV-DN--ADMM)} \quad
\left\{
\begin{aligned}
\min_{y,\, \Delta b^M} \quad
& f(y)
+ \sum_{k \in \mathcal{K}} \sum_{t \in \mathcal{T}}
\mu^{s}_{k,t}
\left(
Mp^{s}_{k,t} - Lp^{s,(r)}_{k,t}
\right) \\
& \quad
+ \sum_{k \in \mathcal{K}} \sum_{t \in \mathcal{T}}
\frac{\rho}{2}
\left\|
Mp^{bc}_{k,t}
- Lp^{bc,(r)}_{k,t}
- u^{bc,(r)}_{k,t}
\right\|^{2} \\[0.5em]
\text{s.t.} \quad
& g(y) = 0, \\
& h(y) \le 0, \\
& \Delta b^{M}_{k,t}
= -\Delta b^{L,(r)}_{k,t},
\quad
\forall k \in \mathcal{K},\ \forall t \in\mathcal{T}.
\end{aligned}
\right.
\end{equation}

In \eqref{LV_ADMM_problem} and \eqref{MV_ADMM_problem}, the upper-level LV-DN subproblems and the lower-level MV-DN problem are reformulated into a structure that can be coordinated through an ADMM-based consensus mechanism combined with Lagrangian relaxation. The dual variables $\mu^{s}_{k,t}$ and $u^{bc,(r)}_{k,t}$ denote the Lagrange multipliers associated with the relaxed consistency constraints for the surplus-energy exchange and the cheap-energy exchange, respectively. The superscript $(r)$ denotes the ADMM iteration index, indicating that the corresponding variables are treated as fixed parameters obtained from the previous iteration when solving each subproblem. The parameter $\rho > 0$ is the ADMM penalty parameter associated with the augmented Lagrangian term, which controls the enforcement of consensus among the shared cheap-energy exchange variables and influences the algorithm's convergence behavior. 

\vspace{6pt}
\begin{algorithm}[H]
    \caption{LDD-ADMM algorithm}\label{algoritmo}
\begin{algorithmic}[1]
\State Initialize $Mp^{bc,(0)}_{k,t}$, $u^{bc,(0)}_{k,t}$, $\mu^{s,(0)}_{k,t}$, tolerance $\varepsilon$, iteration counter $r \gets 0$
 \State \textbf{While }$\text{GAP} > \varepsilon$:
    \State \quad \( r \gets r + 1 \)
    \State \quad \textbf{For} each $k \in \mathcal{K}$:
    \State \quad \quad Solve the $LV_k$ using $Mp^{bc,(r-1)}_{k,t}$, $u^{bc,(r-1)}_{k,t}$ and obtain the current solution $Lp^{bc}_{k,t}$,  $Lp^{s}_{k,t}$,  $\Delta b^{L}_{k,t}$: 
    \State \quad Update the MV parameters: $Lp^{s,(r)}_{k,t} \gets Lp^{s}_{k,t}$, $\enspace Lp_{k,t}^{bc,(r)} \gets Lp^{bc}_{k,t}$, $\enspace \Delta b^{L,(r)}_{k,t} \gets \Delta b^{L}_{k,t}$
    \State \quad Solve the MV and obtain $Mp^{bc,(r)}_{k,t}$ and $Mp^{s,(r)}_{k,t}$
    \State \quad Update parameters $u_{k,t}^{bc,(r+1)}$ and $\mu_t^{s,(r+1)}$:
    \Statex \quad \quad $u_{k,t}^{bc,(r+1)} = u_{k,t}^{bc,(r)} + \tau \, \big(L p_{k,t}^{bc,(r+1)} - M p_{k,t}^{bc,(r+1)}\big)$
    \Statex \quad \quad $\mu_{k,t}^{s,(r+1)} = \mu_{k,t}^{s,(r)} + \eta \big(M p_{k,t}^{s,(r+1)} - L p_{k,t}^{s,(r+1)}\big)$
    \State \quad Update $LV_k$ parameters 
    \State \quad Compute primal and dual residuals:
    \Statex \quad \quad $r^{(r+1)} = \left(\sum_{k,t} \big(L p_{k,t}^{bc,(r+1)} - M p_{k,t}^{bc,(r+1)}\big)^2 \right)^{1/2}$
    \Statex \quad \quad $s^{(r+1)} = \rho \left(\sum_{k,t} \big(M p_{k,t}^{bc,(r+1)} - M p_{k,t}^{bc,(r)}\big)^2 \right)^{1/2}.$
    \State \quad Compute stop criteria $\text{GAP}^{(r+1)} = \max(r^{(r+1)}, s^{(r+1)})$
\State \textbf{End}
\end{algorithmic}
\end{algorithm}
\vspace{6pt}

Starting from Step~4, the algorithm in~\ref{algoritmo} solves in parallel, for each $k \in \mathcal{K}$, the LV subproblem $LV_k$, obtaining the optimal boundary exchange decisions $Lp^{bc}_{i,t}$, $Lp^{s}_{i,t}$, and the net boundary power $\Delta b^L_{i,t}$. These quantities are then communicated to the MV problem and treated as fixed parameters in Step~6. In Step~7, the MV problem is solved using the aggregated LV decisions, yielding the available cheap energy allocation $Mp^{bc,(r)}_{k,t}$ for each LV network. Then, in Step~8, the ADMM coordination parameters are updated. Specifically, $u^{bc}_{k,t}$ corresponds to the scaled dual variable associated with the consensus constraint on cheap energy exchanges, while $\mu_{k,t}^{s}$ is the Lagrange multiplier related to surplus energy balance. These updates enforce consistency between the cheap energy requested by each LV network and the amount supplied by the MV system. Subsequently, in Step~9, the updated values of $Mp^{bc}_{k,t}$ and $u^{bc}_{k,t}$ are passed back to the $LV_k$ subproblems for the next iteration.

In Step~10, the primal residual $r^{(r+1)}$ and the dual residual $s^{(r+1)}$ are computed to quantify, respectively, the violation of the consensus constraints and the variation of the shared variables between successive iterations. Finally, the stopping criterion is evaluated in Step~11 through the gap metric $\text{GAP}^{(r+1)} = \max\!\left(r^{(r+1)}, s^{(r+1)}\right)$, and the algorithm terminates when this value falls below the prescribed tolerance $\varepsilon$.  Finally, to enhance numerical stability and mitigate oscillatory behavior, relaxed dual updates are used in the coordination procedure. In particular, the scaled dual variable associated with the cheap-energy consensus constraint is updated using a damping factor $\tau$, so that the multiplier responds gradually to the mismatch between the cheap energy requested by the LV networks and the amount allocated by the MV network. Similarly, the parameter $\eta$ controls the update of the Lagrange multipliers associated with the surplus-energy exchange constraints \cite{goncalves2017convergence,wei2022tfpnp}. In the implementation, these parameters are set to $\tau=0.5$ and $\eta=0.05$, respectively, thereby reducing excessively aggressive multiplier corrections and improving the numerical behavior of the iterative scheme. Additionally, the ADMM penalty parameter is set to $\rho=15$, while the convergence tolerance is fixed at $\varepsilon=10^{-5}$ and the maximum number of iterations is set to 30.


\section{Case study and computation results}\label{sec:Case_study}
This section presents the case study and computational assessment of the proposed bilevel coordination framework. The test system configuration is first described, followed by a comparison between the exact KKT-based reformulation and the distributed LDD-ADMM algorithm in terms of optimality gap and computational performance. Finally, the structural characteristics of the obtained solutions are analyzed under incremental LV integration scenarios.

\subsection{Case Study}
The proposed bottom-up coordination framework, including the bilevel model and the LDD-ADMM distributed algorithm, was tested using a modified version of the IEEE 33-bus system \cite{9258930} to represent the MV network, and a reduced 206-bus version of the European LV test feeder \cite{PIGEM2026120031} to represent each LV network. These test systems were selected for their widespread use in distribution system studies and for their ability to capture operational constraints at both voltage levels.

The system comprises 33 buses and 32 distribution lines, with peak active and reactive loads of 2.345 MW and 1.310 MVAr, respectively. The MV grid includes 23 conventional load buses, a maximum dispatchable generation capacity of 16 MW, and a utility-scale PV installation with a maximum capacity of 1 MW. A total of six LV networks are connected to the MV system through dedicated interface buses. The main characteristics of the MV network are summarized in Table \ref{tab:MV_data}. Likewise, each LV network corresponds to a 206-bus radial feeder and is connected to a specific MV bus, as indicated by the column labels (Bus 3, Bus 8, Bus 19, Bus 21, Bus 24, and Bus 29). Table \ref{tab:LV_data} reports the characteristics of the six LV networks connected to the MV system. This configuration reflects the physical interface between voltage levels, in which each LV community exchanges power with the MV grid at a single point of common coupling.

\begin{table}[h]
\centering
\caption{Main characteristics of the medium-voltage distribution network.}
\label{tab:MV_data}
\begin{tabular}{lc}
\hline
Parameter & Value \\
\hline
Number of buses & 33 \\
Number of lines & 32 \\
Peak active load (MW) & 2.345 \\
Peak reactive load (MVAr) & 1.310 \\
Maximum conventional generation (MW) & 16 \\
Maximum utility-scale PV capacity (MW) & 1 \\
Number of conventional loads & 23 \\
Number of connected LV networks & 6 \\
\hline
\end{tabular}
\end{table}

\begin{table*}[h]
\centering
\caption{Main characteristics of the low-voltage distribution networks connected to the MV system.}
\label{tab:LV_data}
\begin{tabular}{lcccccc}
\hline
Parameter & Bus 3 & Bus 8  & Bus 19 & Bus 21 & Bus 24 & Bus 29 \\
\hline
Number of buses & 206 & 206 & 206 & 206 & 206 & 206 \\
Number of lines & 205 & 205 & 205 & 205 & 205 & 205 \\
Number of users & 23 & 50 & 23 & 22 & 45 & 23 \\
Number of PV units & 9 & 22 & 12 & 11 & 19 & 9 \\
Number of BESS units & 2 & 4 & 3 & 3 & 5 & 2 \\
Peak active load (MW) & 0.078 & 0.170 & 0.078 & 0.075 & 0.153 & 0.078 \\
Peak reactive load (MVAr) & 0.023 & 0.051 & 0.023 & 0.022 & 0.046 & 0.023 \\
Total daily energy demand (MWh) & 0.792 & 1.716 & 0.803 & 0.749 & 1.518 & 0.792 \\
Installed PV capacity (MW) & 0.051 & 0.110 & 0.060 & 0.055 & 0.095 & 0.051 \\
Installed BESS capacity (MWh) & 0.010 & 0.020 & 0.015 & 0.015 & 0.025 & 0.010 \\
PV penetration (\%) & 65.2 & 64.7 & 76.7 & 73.5 & 62.1 & 65.2 \\
BESS energy capacity ratio (\%) & 12.8 & 11.8 & 19.2 & 20.1 & 16.3 & 12.8 \\
\hline
\end{tabular}
\end{table*}

To introduce heterogeneity of the analysis, each LV network was parameterized with different levels of demand, PV penetration, and BESS capacity. Peak active demand ranges from 0.075 MW to 0.170 MW, while peak reactive demand varies accordingly. Installed PV capacity and BESS capacity were also differentiated across LV networks, resulting in PV penetration levels of approximately 62\% to 77\% and varying energy storage capacity ratios. This diversity allows the coordination mechanism to be evaluated under non-uniform operating conditions, reflecting differences in consumption patterns and DER adoption across communities. All simulations were conducted over a 24-hour horizon with a time resolution of 1 hour. Distinct demand and generation profiles were assigned to each network to capture temporal variability. The detailed network data, parameter settings, and load and generation profiles are provided in the supplementary material to ensure reproducibility.

To analyze the structural impact of LV integration on the MV system, an incremental integration exercise was performed. Starting from a single connected LV network, additional LVs were progressively incorporated into the MV system until all six communities were connected. This stepwise expansion allows the marginal effect of each new LV to be isolated with respect to system operation, resource allocation, and computational performance. Table \ref{tab:LV_energy_share} reports the relative demand share of each LV network under incremental integration. As additional communities are connected, the aggregate LV demand increases, and the relative weight of each LV adjusts accordingly. These demand shares play a central role in the redistribution of the capacity-limited PV resource at the MV level, as analyzed in subsequent subsections.

\begin{table}[h]
\centering
\caption{Relative energy share (\%) of each LV-DN under incremental integration into the MV network.}
\label{tab:LV_energy_share}
\begin{tabular}{lcccccc}
\hline
LV-DN & 1 LV & 2 LVs & 3 LVs & 4 LVs & 5 LVs & 6 LVs \\
\hline
Bus 3  & 100.0 & 31.6 & 23.9 & 19.5 & 14.2 & 12.4 \\
Bus 8  & 0.0   & 68.4 & 51.8 & 42.3 & 30.8 & 26.9 \\
Bus 19 & 0.0   & 0.0  & 24.3 & 19.8 & 14.4 & 12.6 \\
Bus 21 & 0.0   & 0.0  & 0.0  & 18.4 & 13.4 & 11.8 \\
Bus 24 & 0.0   & 0.0  & 0.0  & 0.0  & 27.2 & 23.8 \\
Bus 29 & 0.0   & 0.0  & 0.0  & 0.0  & 0.0  & 12.4 \\
\hline
\end{tabular}
\end{table}

The models were implemented in Python using Pyomo and solved with Gurobi. All experiments were executed on a workstation equipped with an Intel Core i9-14900HX processor (2.20 GHz) and 64 GB of RAM under a 64-bit operating system.

\subsection{Distributed coordination algorithm performance}

Table~\ref{Table_Comparison_KKT_ADMM} compares the solution obtained through the exact KKT reformulation, where complementarity conditions are enforced using SOS1 constraints, and the proposed LDD-ADMM algorithm for an increasing number of LV networks connected to the MV system (up to six LVs). The table reports the objective function value, the total number of variables (including binary variables), and the computational time for both approaches. For the LDD-ADMM algorithm, the optimality gap, primal and dual residuals, number of iterations, and the relative difference with respect to the KKT solution are also provided. Additionally, the last three columns report the accuracy of the second-order cone (SOC) relaxation used to model the DN constraints, computed as in \cite{9535394}. The LVs column corresponds to the maximum SOC relaxation violation across all LV networks. For the MV network, MV* denotes the SOC relaxation violation computed ex post, considering network losses, and is used as the reference value, whereas MV denotes the corresponding SOC relaxation violation obtained using the lossless MV approximation adopted in the optimization model. In all tested cases, the reported violations remain below $2.0\mathrm{E}{-06}$ for the LV networks and below $7.7\mathrm{E}{-04}$ for the MV network, suggesting that the SOC relaxation provides a sufficiently accurate representation of the underlying network constraints.

\begin{table}[h]
\centering \small
\setlength{\tabcolsep}{3pt}
\begin{tabular}{c|cccc|ccccccc|ccc}
\hline
 & \multicolumn{4}{c|}{\textbf{KKT}} 
 & \multicolumn{7}{c|}{\textbf{LDD-ADMM}}
 & \multicolumn{3}{c}{\textbf{SOC relaxation}} \\ 
\hline
\textbf{LVs} 
& \textbf{Obj} 
& \textbf{Vars} 
& \textbf{Bin} 
& \textbf{t (s)}
& \textbf{Obj} 
& \textbf{Gap} 
& \textbf{$r$} 
& \textbf{$s$} 
& \textbf{t (s)} 
& \textbf{Iter} 
& \textbf{Rel. Gap}
& \textbf{LVs}
& \textbf{MV*}
& \textbf{MV} \\ 
\hline
1 & 11.1 & 49488  & 5736  & 1.5  & 11.1 & 3.0E-06 & 3.0E-06 & 4.0E-07 & 37.9  & 17 & 3.7E-04 & 3.9E-07 & 3.3E-06 & 5.8E-05 \\
2 & 35.6 & 82272  & 7056  & 2.9  & 35.6 & 1.6E-06 & 1.6E-06 & 3.0E-07 & 100.6 & 26 & 7.7E-05 & 2.0E-06 & 6.1E-05 & 6.3E-05 \\
3 & 46.6 & 115056 & 8376  & 4.1  & 46.6 & 9.1E-06 & 9.1E-06 & 9.1E-06 & 85.8  & 18 & 1.0E-04 & 2.0E-06 & 9.4E-05 & 9.5E-05 \\
4 & 57.0 & 147840 & 9696  & 6.7  & 57.0 & 1.6E-05 & 1.6E-05 & 6.8E-06 & 161.5 & 30 & 3.5E-04 & 2.0E-06 & 4.0E-04 & 4.1E-04 \\
5 & 79.4 & 180624 & 11016 & 7.4  & 79.4 & 6.0E-06 & 1.8E-06 & 6.0E-06 & 170.1 & 25 & 6.7E-05 & 2.0E-06 & 6.3E-04 & 6.4E-04 \\
6 & 91.1 & 213408 & 12336 & 11.1 & 91.1 & 8.2E-06 & 2.2E-06 & 8.2E-06 & 172.1 & 23 & 6.1E-05 & 2.0E-06 & 7.6E-04 & 7.7E-04 \\
\hline
\end{tabular}
\caption{Comparison between the KKT reformulation and the proposed LDD-ADMM algorithm for an increasing number of LV networks, including SOC relaxation accuracy metrics.}\label{Table_Comparison_KKT_ADMM}
\end{table}

The most relevant result concerns the solution quality delivered by the LDD-ADMM algorithm. Despite being theoretically inexact due to the dualization of coupling constraints and the iterative coordination process, the algorithm yields solutions close to those obtained from the exact KKT reformulation, with deviations remaining within a numerical tolerance. The reported optimality gaps and primal/dual residuals remain on the order of $10^{-6}$ across all tested cases, indicating strong numerical convergence. Furthermore, the relative deviation of the objective function value remains below $0.04\%$ across all cases. These results show that the proposed LDD-ADMM algorithm preserves solution accuracy while relaxing the centralized structure of the exact reformulation. In this regard, it is important to emphasize that the LDD-ADMM algorithm is not designed to improve optimality but to enable distributed coordination while enhancing LV-level privacy. In contrast, the KKT reformulation yields a fully centralized problem in which all LV information must be explicitly revealed. The near-equivalence of the objective values therefore shows that decentralization is achieved without sacrificing economic optimality.

The computational cost of decentralization is reflected in the solution times. On average, the LDD-ADMM algorithm requires approximately one order of magnitude more computational time than the KKT reformulation. Nevertheless, all solution times remain within a day-ahead operational timeframe. Hence, the additional computational burden can be interpreted as the cost of enabling distributed coordination and enhancing privacy rather than a loss in solution quality. However, the number of ADMM iterations does not increase monotonically with the number of LVs. The algorithm converges in fewer than 30 iterations in all cases, which is acceptable given the achieved accuracy. This suggests that the convergence behavior is stable with respect to system size and that the coordination mechanism scales appropriately under the tested conditions.

\begin{figure}[h]
\centering
\includegraphics[width=6.5in]{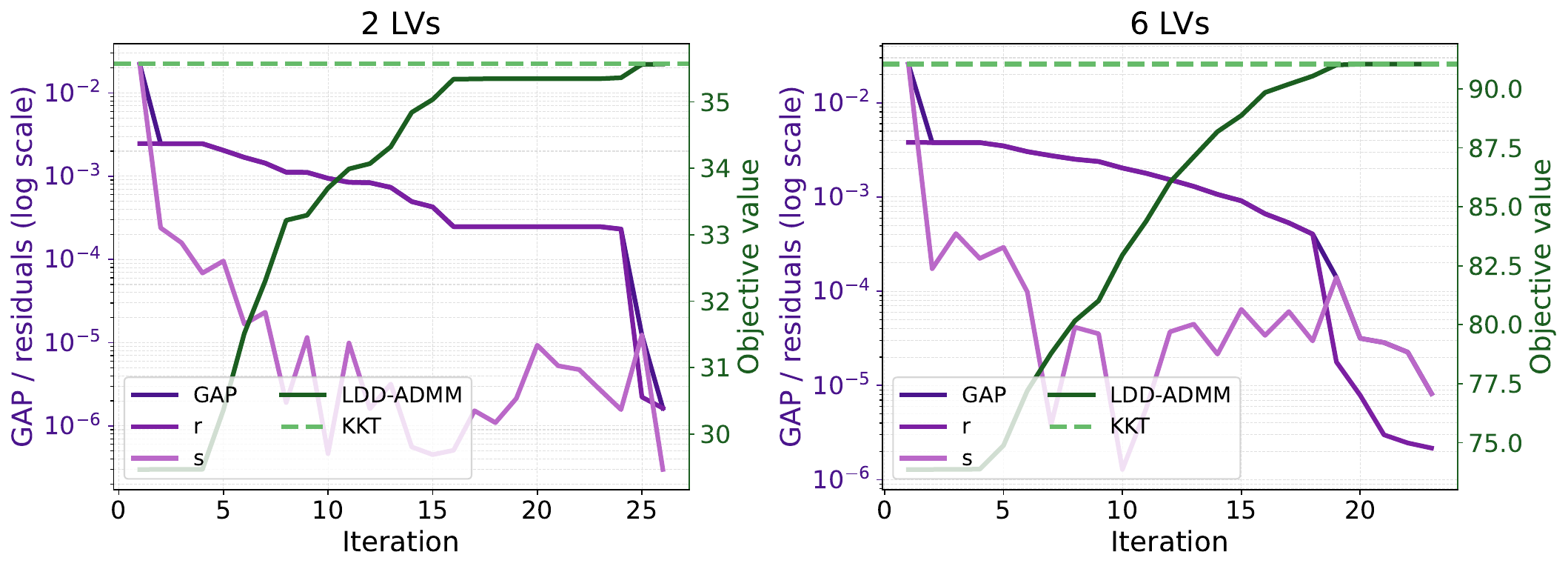}
\caption{Convergence behavior of the proposed LDD-ADMM algorithm for two and six LV communities.}
\label{Fig_Convergence}
\end{figure}

Figure~\ref{Fig_Convergence} illustrates the convergence behavior of the proposed LDD-ADMM algorithm for two and six LV communities. In both instances, the optimality gap and primal/dual residuals decrease from the order of $10^{-2}$ to below $10^{-5}$ within fewer than 30 iterations, confirming stable numerical convergence. Two convergence phases can be distinguished: (i) an initial adjustment phase (approximately the first 10 iterations), characterized by a rapid reduction of the gap and significant updates in the objective value as coordination among LV and MV subproblems is established; and (ii) a refinement phase, in which the residuals progressively decrease and the objective value approaches the reference KKT solution through minor incremental corrections. Importantly, the convergence pattern remains qualitatively unchanged as the number of LV communities increases from two to six. The number of iterations required to reach the prescribed tolerance does not increase, indicating that the coordination mechanism scales without degradation in convergence behavior. Moreover, the objective value evolves smoothly and monotonically toward the centralized KKT benchmark, with no observable oscillations or instability.

\subsection{Dispatch Equivalence and Resource Allocation Differences}

In the previous subsection, the performance of the proposed algorithm was evaluated with respect to optimality and convergence to the global solution. In this subsection, we analyze differences in the solution structure, i.e., the specific combination of decision variables, rather than focusing exclusively on the optimal objective value.

\begin{figure}[h]
\centering
\includegraphics[width=0.99\textwidth]{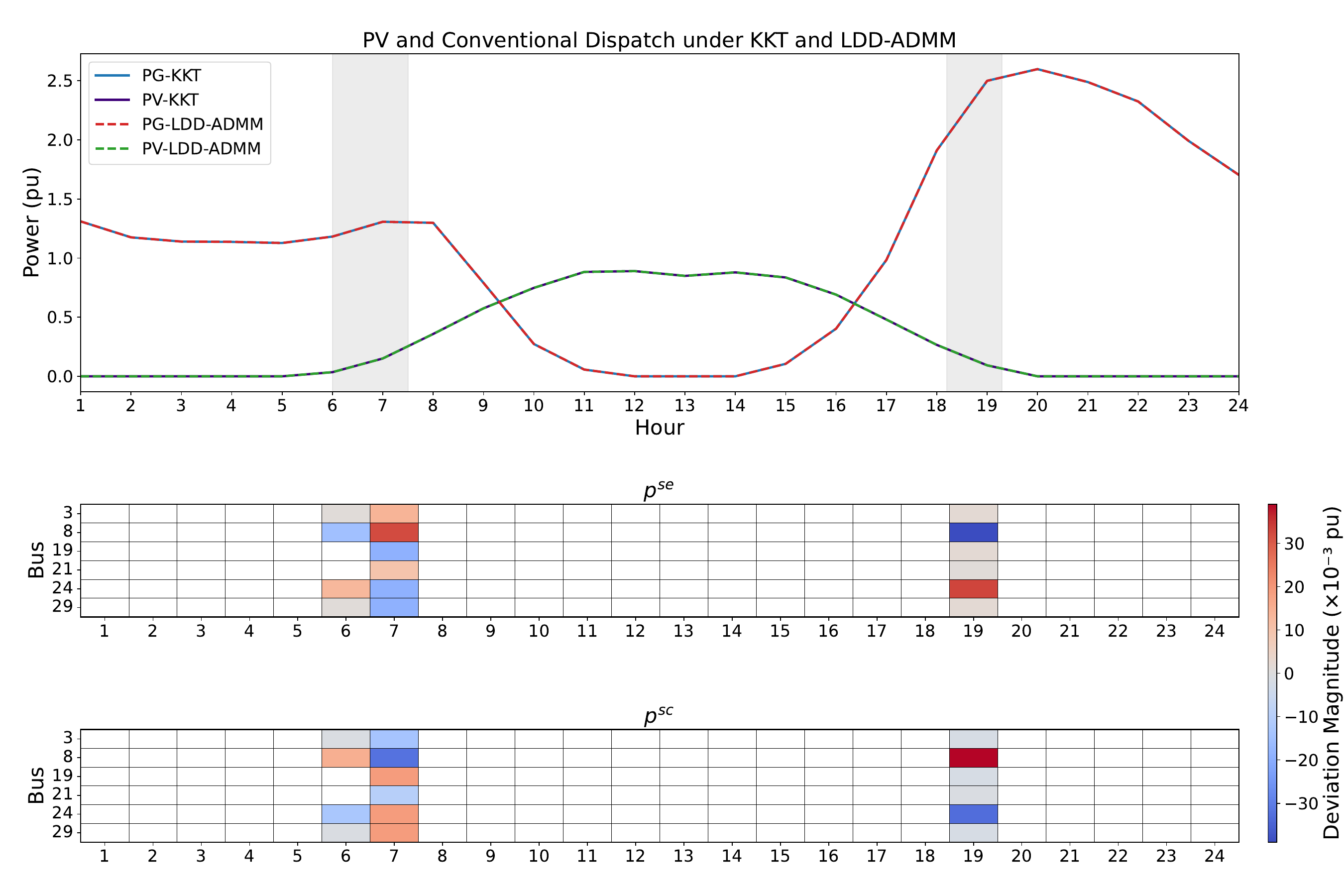}
\caption{Dispatch equivalence and structural deviations between KKT reformulation and LDD-ADMM solution.}
\label{Fig_Comparison}
\end{figure}

Figure~\ref{Fig_Comparison} summarizes the comparison between the KKT reformulation and the proposed LDD-ADMM algorithm. The upper panel compares the dispatch of conventional (expensive) and PV (cheap) generation obtained using the KKT reformulation and the proposed LDD-ADMM algorithm. The lower panels illustrate how this generation is allocated among the LV communities by reporting the differences in the exchanged quantities $p^{se}$ and $p^{sc}$ between both solution approaches.

Focusing first on the upper panel, the dispatch curves obtained with both approaches are perfectly superimposed for both conventional and PV generation. This behavior is theoretically expected. The follower’s problem corresponds to a convex cost-minimization problem. Under convexity and regularity conditions (e.g., Slater’s condition), the KKT conditions are both necessary and sufficient for optimality. Therefore, any solution obtained via the KKT reformulation must coincide with the primal optimum of the follower’s problem. Since the distributed LDD-ADMM algorithm solves the same convex subproblem through dual decomposition and consensus enforcement, both approaches converge to the same optimal dispatch. The superposition of the curves thus validates the correctness of the distributed implementation and confirms that the follower’s minimum-cost solution is preserved.

The lower panels of Figure~\ref{Fig_Comparison} present heatmaps depicting the differences in $p^{se}$ and $p^{sc}$, where each value corresponds to the deviation of the LDD-ADMM solution from its KKT counterpart. As a result, the white cells indicate zero difference. The variable $p^{b}$, together with the net boundary exchange $\Delta b$, is not reported because all differences are exactly zero across all buses and time periods. This exact agreement at the aggregate level implies that the quantities governing the coupling between LV and MV networks are identical under both solution approaches, ensuring that the leader–follower interaction and the resulting equilibrium remain unchanged.

When the aggregate energy exchange is disaggregated, differences emerge only in variables related to the allocation of cheap and expensive energy among the LV communities. As shown in the lower panels of Figure~\ref{Fig_Comparison}, these deviations are concentrated in two narrow time intervals corresponding to the beginning and the end of PV production, when the availability of the cheap resource becomes scarce. Outside these periods, both approaches produce identical allocations across all communities and hours. The magnitude of the observed deviations remains on the order of $10^{-3}$ and exhibits a symmetric pattern between $p^{se}$ and $p^{sc}$. For instance, a value of $-0.015$ in $p^{se}$ corresponds to $+0.015$ in $p^{sc}$. One possible explanation is that the proposed LDD-ADMM algorithm provides an approximate distributed solution of the exact KKT-SOS1 reformulation. Consequently, small residual discrepancies may remain in some allocation variables when the stopping criterion is reached. The fact that these deviations occur only during periods when the cheap-resource constraint is active suggests that the consensus process is most sensitive under these operating conditions. Nevertheless, their magnitude remains negligible, consistent with the reported optimality gap of 0.04\%, and does not affect the aggregate energy exchanges or the economic outcomes of the coordination process.


\subsection{Operational Impact of the Stackelberg Coordination Structure}
To assess the operational implications of the proposed Stackelberg framework, the MV dispatch induced by the bilevel interaction was compared with the dispatch obtained when the MV system operates independently. Since the upper-level objective accounts only for LV operating costs, MV operating costs are not directly represented in the bilevel objective. Therefore, after solving the bilevel problem, the resulting imports and exports were imposed as fixed exchange conditions in the standalone MV optimization model, allowing the follower's operating cost to be recovered. Under this imposed operational allocation, the MV-level objective function reached 56.22 \$, which corresponds to the operating cost of the MV system under the bilevel solution.

Subsequently, to assess the economic outcome for the MV system under autonomous operation, the standalone MV optimization problem was solved by enforcing only the aggregated net exchange profile $\Delta b$ obtained from the bilevel solution. In contrast to the previous analysis, the individual import and export allocations were not fixed, allowing the MV system to optimally allocate its available resources according to its own objective function. The resulting MV operating cost was 50.07 \$, lower than that obtained under the fully prescribed exchange schedule. This difference reveals that the operational response induced by the bilevel interaction is the cost of allocating generation resources that the MV level would naturally select under an isolated congestion- and efficiency-oriented operational criterion.

\begin{figure}[htpb]
\centering
\includegraphics[scale=0.50]{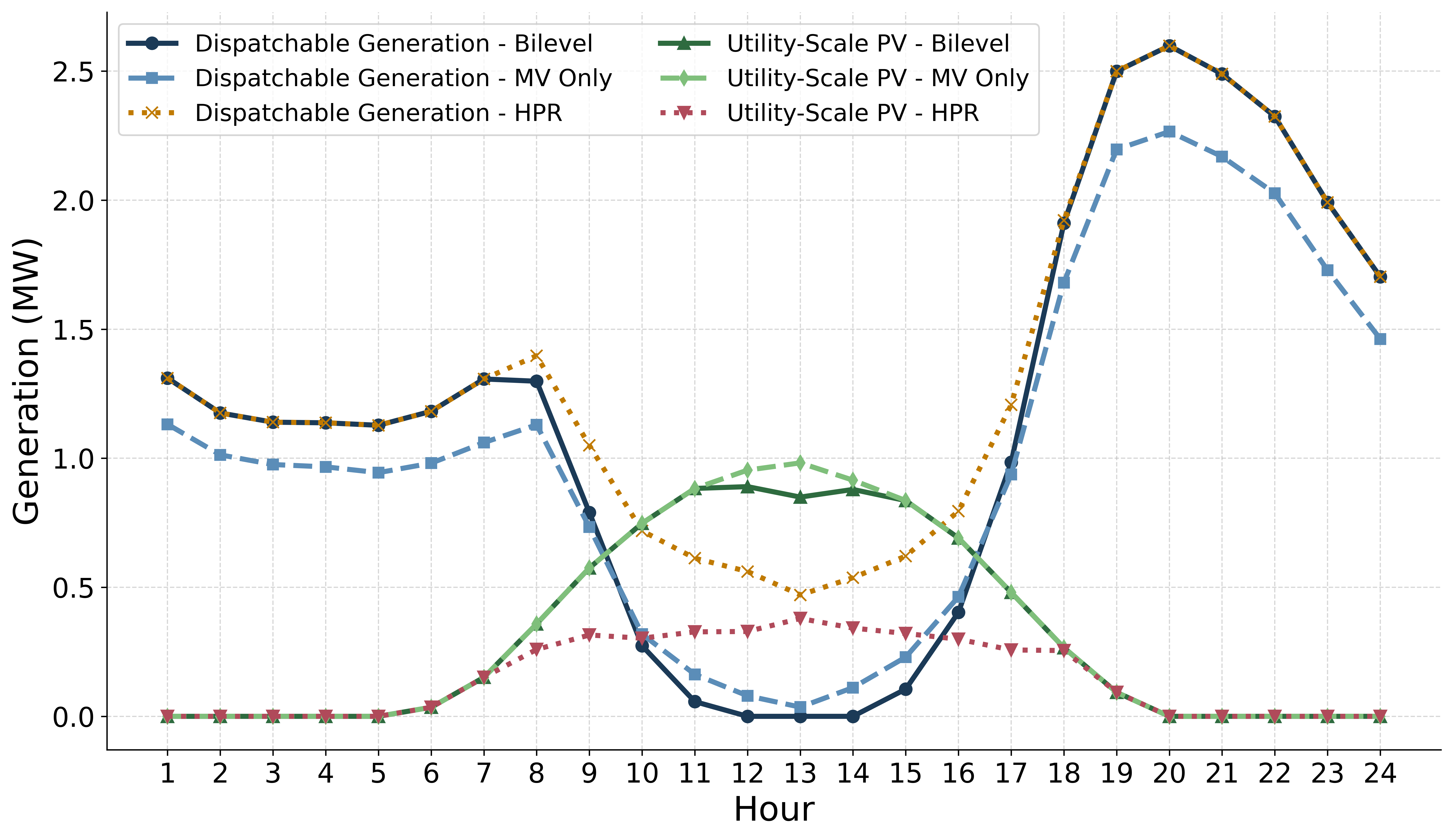}
\caption{Generation Scheduling Comparison.}
\label{Generation_Diference_Fig}
\end{figure}

To support this distinction, Fig. \ref{Generation_Diference_Fig} compares the generation schedules obtained from the bilevel solution and the isolated MV optimization. The results show that when the MV system is allowed to allocate its available resources freely, dispatch relies more on utility-scale PV generation during periods of high renewable availability, reducing the utilization of dispatchable conventional generation. In contrast, under the bilevel solution, the dispatchable generation remains systematically higher, particularly during the transition periods between renewable-dominated and thermal-dominated operation. This behavior indicates that the LV-level exchange requirements partially restrict the MV system's operational flexibility, preventing the follower from reaching its preferred generation allocation. From a network operations perspective, this effect is also influenced by the electrical topology of the MV system and the spatial distribution of LV communities across the network. Due to the resistive nature of distribution systems and their relatively R/X ratios, the allocation of cheap and expensive resources cannot be interpreted exclusively from an energy-balance perspective, since the geographical location of injections and withdrawals affects voltage profiles, line loadings, and network losses. Consequently, the MV system may require additional dispatchable generation to maintain secure operating conditions even when low-cost renewable resources remain available elsewhere in the network.

To further examine the role of the lower-level problem, a High-Point Relaxation (HPR) of the bilevel model was evaluated. The HPR was obtained by imposing the follower primal feasibility constraints into the upper level and removing the lower-level objective function, thereby yielding a single-level problem. Under this relaxation, the MV network constraints continue to guarantee the solution's physical and operational feasibility, but no longer enforce the cost-minimizing dispatch that characterizes the follower response in the original bilevel formulation. 

Although the HPR yields the same upper-level objective value as the original bilevel formulation, the resulting MV dispatch differs significantly. Fig. \ref{Generation_Diference_Fig} illustrates this effect. Under the bilevel formulation, utility-scale PV generation is exploited whenever available, while conventional generation is used only to satisfy the remaining demand. In contrast, under the HPR, only a portion of the available PV generation is utilized, with the remaining MV demand being supplied by dispatchable generation. This behavior arises because exchanges between the MV network and the LV communities are directly represented in the upper-level objective via exchange prices, whereas the conventional demand connected at the MV level does not appear in the leader's objective function. Consequently, once the follower objective is removed, the model no longer contains an operational criterion that promotes the efficient allocation of PV resources across all MV demands, allowing conventional generation to replace renewable generation even when the latter remains available. 

This result highlights the role of the bilevel approach within the proposed coordination framework. As DNs evolve towards increasingly heterogeneous environments, where energy communities, conventional consumers, aggregators, and other active participants coexist, the network operator's role extends beyond enforcing network feasibility alone. In such settings, an explicit operational criterion is required to allocate shared generation resources across competing uses while preserving an efficient system-wide dispatch. Otherwise, a single-level formulation based solely on network feasibility may produce physically admissible solutions that do not reflect the resource allocation that would emerge under coordinated operation. This effect is clearly observed in the present case study, where the HPR remains feasible but fails to reproduce the dispatch obtained when utility-scale PV and conventional generation must simultaneously support both LV-community exchanges and conventional MV demand efficiently. Therefore, the HPR results suggest that the bilevel structure becomes relevant as active participants proliferate within DNs and resource allocation decisions become intrinsically coupled across multiple user groups.

\begin{figure}[htpb]
\centering
\includegraphics[scale=0.42]{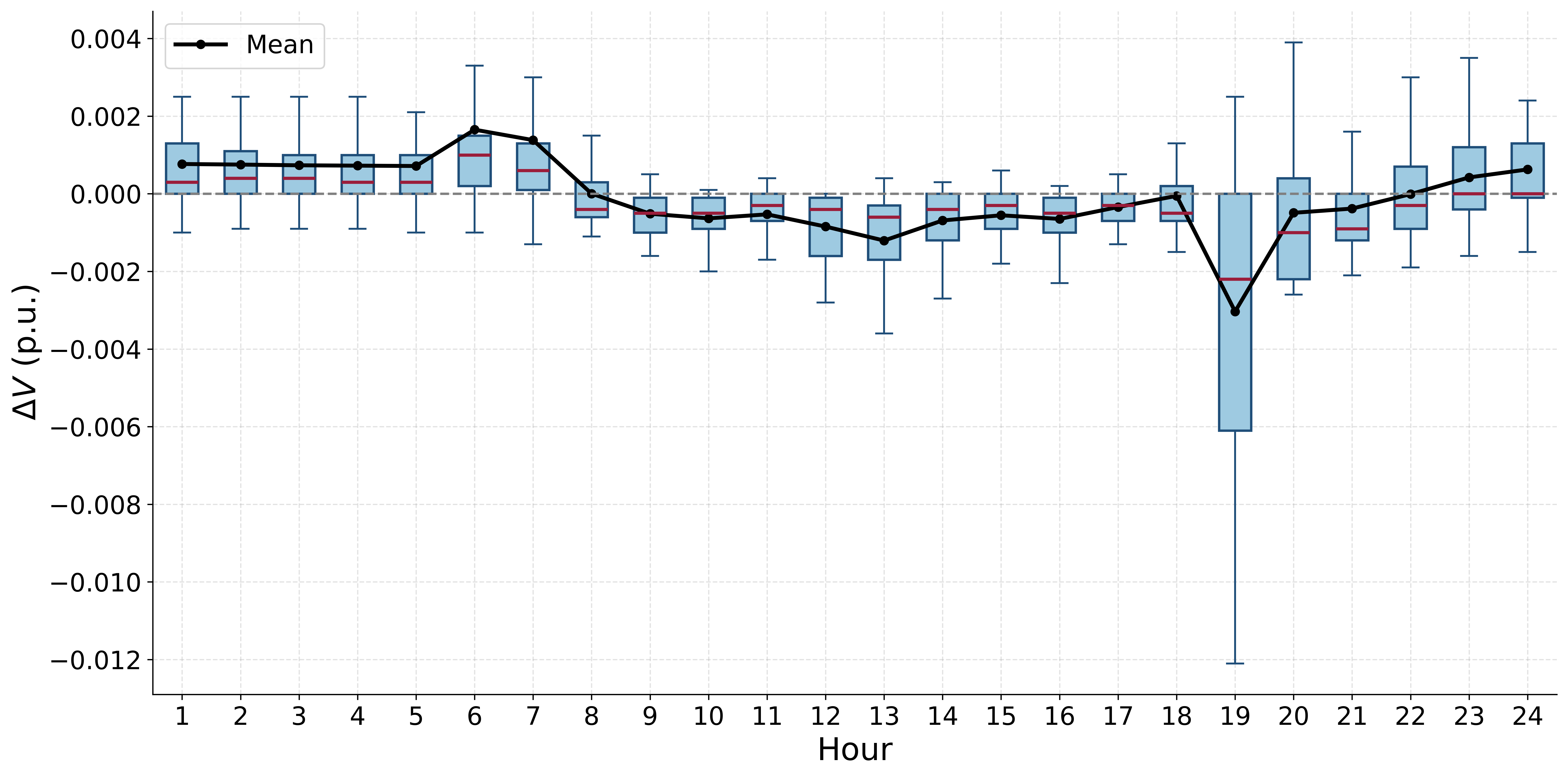}
\caption{Voltage Difference Distribution (Bilevel - MV Only)}
\label{voltage_difference_Fig}
\end{figure}

In the above context, the operational consequences of this interaction can be analyzed in the voltage behavior shown in Fig. \ref{voltage_difference_Fig}, which presents the difference between the voltage profiles obtained under the bilevel solution and those resulting from the isolated MV optimization. Positive values indicate operating points at which the bilevel solution produces higher voltages, while negative values indicate the opposite. Although the voltage deviations remain relatively small, their temporal structure reveals a clear operational pattern associated with the generation schedules observed in Fig. \ref{Generation_Diference_Fig}. During periods of high PV utilization, the bilevel allocation tends to produce slightly higher voltage levels in several network locations due to the redistribution of renewable injections across the MV grid. Conversely, during evening hours, when the system becomes increasingly dependent on dispatchable generation, the bilevel interaction induces localized voltage reductions relative to the isolated MV solution. These patterns indicate that the Stackelberg interaction does not merely alter energy exchanges at the aggregate level but also modifies the spatial distribution of power flows throughout the network.

\begin{figure}[htpb]
\centering
\includegraphics[scale=0.4]{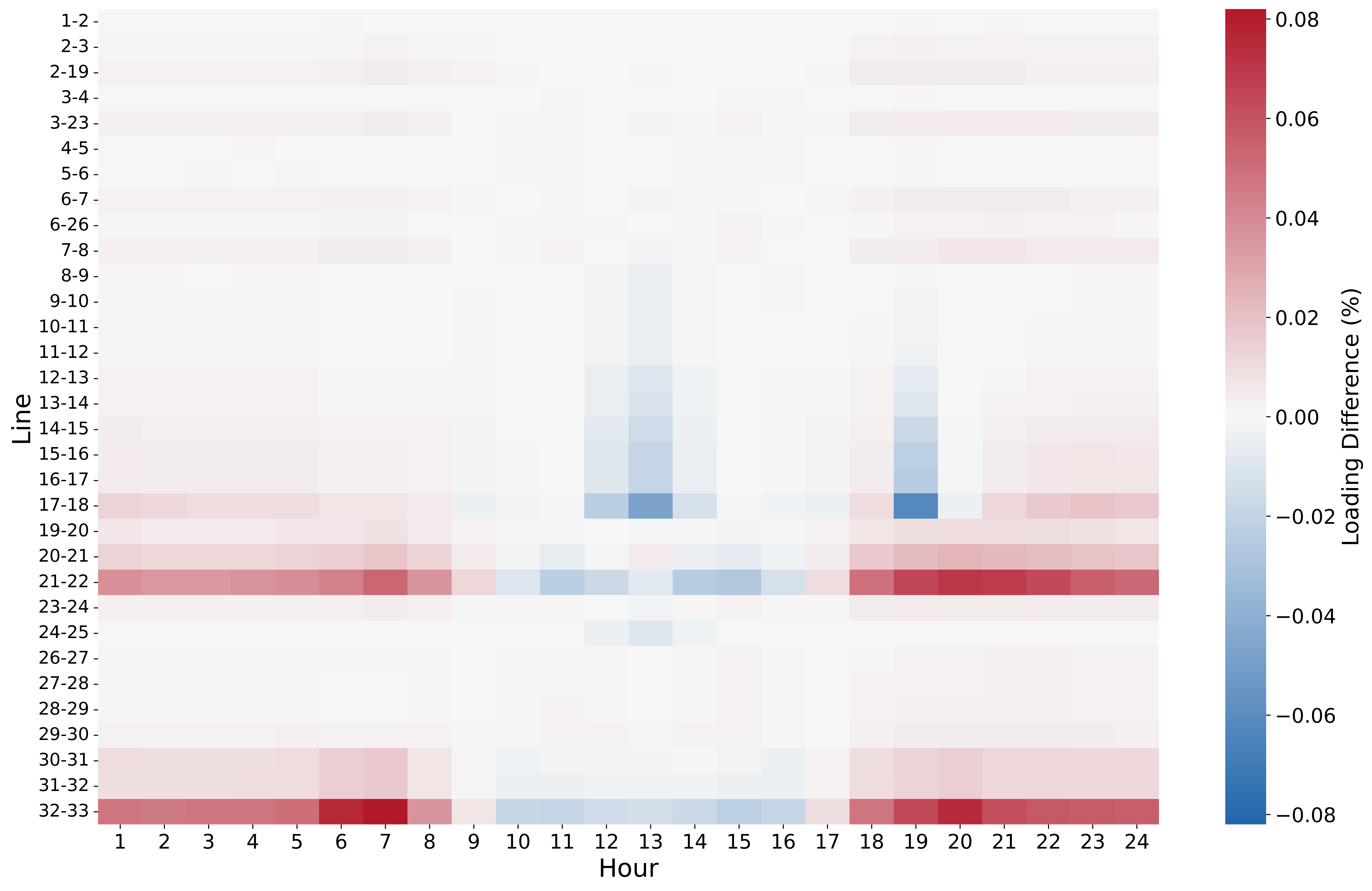}
\caption{Line Loading Difference (Bilevel - MV Only)}
\label{loading_difference_Fig}
\end{figure}

Finally, Fig. \ref{loading_difference_Fig} further confirms this effect by showing the differences in line loading between the two operating conditions. Although most branches exhibit relatively small variations, the bilevel solution generates localized increases in loading in specific network sections and during particular operating periods. These increments are not uniformly distributed across the network, but instead appear in branches electrically associated with the redistribution of renewable and dispatchable resources required to satisfy the exchange structure imposed by the LV level. This behavior supports the view that the follower problem is not merely a feasibility-verification layer. Rather, the MV system adapts its operational allocation to accommodate the preferences of the upper-level coordinator, even when those decisions differ from the allocation that would minimize congestion and redispatch requirements under a purely system-oriented operation criterion. Overall, these results empirically reinforce the strategic interpretation previously discussed in the paper, showing that the proposed bilevel framework induces measurable operational effects on the MV network and captures the intrinsic interaction between LV autonomy objectives and MV operational efficiency.

\subsection{Impact of Incremental LV Integration on Cheap Resource Allocation}

This section evaluates how the allocation of the cheap resource (utility-scale PV generation at the MV level) evolves as additional LV communities are progressively connected to the MV network. The analysis quantifies how the progressive sharing of a capacity-limited PV resource reshapes its distribution across LV networks and modifies the overall MV energy balance.

\begin{figure}[h]
\centering
\includegraphics[width=0.9\textwidth]{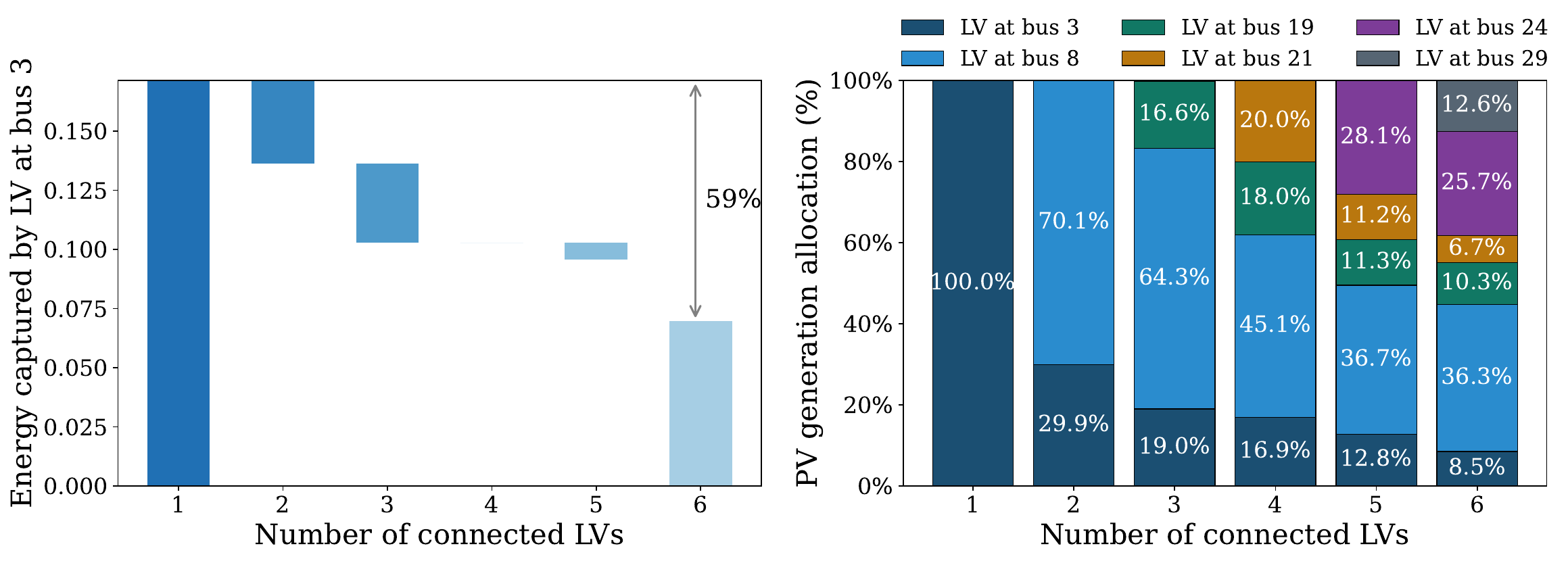}
\caption{Redistribution of capacity-limited PV generation under incremental LV integration.}
\label{combined_figure}
\end{figure}

Figure \ref{combined_figure} summarizes the main effect of this incremental integration. The left panel shows the evolution of the PV energy captured by the LV connected at bus 3 as new LV communities are incorporated. In the single-LV case, this community captures the corresponding available PV allocation. However, as additional LVs are connected, the PV assigned to bus 3 steadily decreases. The waterfall representation clearly quantifies this effect: when all six LV networks are connected, the PV allocation to bus 3 is reduced by approximately 59\% relative to the initial scenario. This reduction reflects the redistribution of a fixed, low-cost resource among an increasing number of competing LV communities. In practical terms, the early-connected LV experiences a dilution of individual benefit as access to the limited PV generation becomes shared.

The right panel of Figure \ref{combined_figure} provides further insight into how this redistribution takes place. The stacked bars show the percentage allocation of MV PV generation among the connected LVs at each integration stage. The distribution is clearly non-uniform. For instance, when two LVs are connected, bus 8 receives 70.1\% while bus 3 retains 29.9\%. With three LVs, the shares become 64.3\%, 19.0\%, and 16.6\%, and as more communities are incorporated, the allocation progressively spreads across them. In this regard, the observed distribution closely follows the relative demand shares of the LV networks reported in Table \ref{tab:LV_energy_share}. Larger-demand LVs (e.g., buses 8 and 24) systematically receive higher percentages of cheap generation, while smaller systems (e.g., buses 3 and 29) capture lower shares. This indicates that, although no explicit proportional allocation rule is imposed in the model, the MV feasibility problem implicitly reallocates the limited PV resource in a manner broadly consistent with demand size. The allocation therefore emerges endogenously from the optimal dispatch, rather than from an exogenous fairness constraint.

Table \ref{tab:energy_balance} complements this graphical analysis by reporting the aggregated power balance and energy composition as LVs are incrementally added. Several structural effects can be identified. First, the contribution of LV-exported surplus energy to the MV system increases from 0.3\% in the single-LV case to 2.3\% when six LVs are connected. Although still modest in magnitude, this increasing share indicates that LV communities progressively contribute more surplus generation to the MV layer. Operationally, this implies that during surplus periods, the corresponding MV nodes transition from being net consumers to temporary net generators, and the injected energy is reallocated to other traditional MV loads. This mechanism effectively supports internal redistribution of power flows and can alleviate localized congestion in the MV grid.

\begin{table*}[h]
\centering
\caption{Aggregated energy composition under incremental LV integration.}
\label{tab:energy_balance}
\setlength{\tabcolsep}{4pt}
\begin{tabular}{lcccccccccccc}
\hline
 & \multicolumn{2}{c}{1 LV} 
 & \multicolumn{2}{c}{2 LVs}
 & \multicolumn{2}{c}{3 LVs}
 & \multicolumn{2}{c}{4 LVs}
 & \multicolumn{2}{c}{5 LVs}
 & \multicolumn{2}{c}{6 LVs} \\
\cline{2-13}
 & MWh & \% & MWh & \% & MWh & \% & MWh & \% & MWh & \% & MWh & \% \\
\hline
\textbf{Total generation} 
& 33.0 & 100 & 34.1 & 100 & 34.5 & 100 & 35.0 & 100 & 35.9 & 100 & 36.4 & 100 \\

Conventional generation
& 25.0 & 75.6 & 25.8 & 75.8 & 26.2 & 75.8 & 26.5 & 75.8 & 27.4 & 76.2 & 27.8 & 76.4 \\

Utility-scale PV
& 8.0 & 24.1 & 7.9 & 23.3 & 7.9 & 22.9 & 7.9 & 22.5 & 7.8 & 21.7 & 7.7 & 21.3 \\

Power purchased from LVs
& 0.1 & 0.3 & 0.3 & 0.9 & 0.4 & 1.3 & 0.6 & 1.7 & 0.7 & 2.1 & 0.8 & 2.3 \\

\hline
\textbf{Total load}
& 33.0 & 100 & 34.1 & 100 & 34.5 & 100 & 35.0 & 100 & 35.9 & 100 & 36.4 & 100 \\

MV load 
& 32.6 & 98.5 & 32.6 & 95.5 & 32.6 & 94.3 & 32.6 & 93.1 & 32.6 & 90.7 & 32.6 & 89.5 \\

Aggregated LV load 
& 0.5 & 1.5 & 1.5 & 4.5 & 2.0 & 5.7 & 2.4 & 6.9 & 3.3 & 9.3 & 3.8 & 10.5 \\

\hline
\textbf{LV imports}
& 0.5 & 61.0 & 1.5 & 60.7 & 2.0 & 59.9 & 2.4 & 59.6 & 3.3 & 59.9 & 3.8 & 60.0 \\

Expensive energy $p^{se}$
& 0.3 & 39.4 & 1.1 & 42.5 & 1.4 & 43.6 & 1.8 & 44.6 & 2.6 & 46.5 & 3.0 & 47.2 \\

Cheap energy $p^{sc}$
& 0.2 & 21.6 & 0.5 & 18.2 & 0.5 & 16.3 & 0.6 & 15.0 & 0.7 & 13.4 & 0.8 & 12.8 \\

\hline
\end{tabular}
\end{table*}

Second, utility-scale PV generation remains nearly constant across scenarios (around 7.7–8.0 MW, representing roughly 21–24\% of total generation), confirming that this resource is fully exploited in all configurations. As additional LV demand is connected (aggregated LV load grows from 1.5\% to 10.5\% of total load), the incremental demand is predominantly covered by conventional generation, whose share increases slightly from 75.6\% to 76.4\%. This confirms that the PV resource is capacity-limited and that system expansion primarily triggers higher utilization of dispatchable generation.

Third, the share of cheap energy allocated to LV imports evolves significantly. While the absolute amount of LV imports increases (from 0.5 MW to 3.8 MW), the proportion corresponding to cheap energy decreases from 21.6\% to 12.8\%. This reflects the progressive pressure from the availability of cheap energy per LV as additional communities are connected. However, when considering the system as a whole, a growing fraction of the available PV generation is progressively reassigned from traditional MV loads toward LV communities. In other words, as more LVs are connected, the MV optimization reallocates part of the previously consumed cheap generation from conventional MV demand to serve LV demand.

This behavior indicates that the MV feasibility layer dynamically reshapes resource allocation to maintain cost-efficient operation under structural demand expansion. The results show that the limited, low-cost resource is not distributed uniformly across LV communities, but rather in proportion to their relative demand weights. As additional LVs are incorporated, the MV optimization reallocates capacity-limited PV generation to match the aggregated demand structure. This demand-driven pattern emerges endogenously from the feasibility-constrained dispatch, without requiring any predefined allocation rule.

\section{Conclusions} \label{sec:Conclusions}

This paper presented a bottom-up coordination framework for the hierarchical interaction between multiple LV energy communities and an MV DN. In the proposed scheme, the LV communities collectively act as leaders within a Stackelberg structure, determining discrete DER schedules and boundary energy exchanges, while the MV network operates as the follower, ensuring efficient system feasibility. The coordination mechanism was formulated as a bilevel optimization problem, in which all discrete decisions are confined to the upper level, and the lower-level MV problem remains convex and continuous, allowing an exact single-level reformulation via the KKT conditions. In addition, an algorithm based on Lagrangian relaxation and ADMM was developed, enabling the decomposition of LV subproblems for parallel execution while enhancing the privacy of individual communities and limiting information exchange to boundary quantities.

The computational results show that the distributed algorithm converges to solutions numerically equivalent to the exact KKT reformulation, achieving convergence in a few iterations with optimality gaps below 0.04\%. Moreover, the MV-level dispatch, including both conventional and PV generation schedules, is identical in both approaches, confirming that the distributed scheme preserves the follower’s optimal response. The proposed framework performs acceptably even for large-scale instances involving multiple LV networks (6 feeders, each with 206 buses) coordinated with a 33-bus MV system. Although the distributed scheme entails additional computational time due to the iterative coordination mechanism, solution times remain compatible with day-ahead operational requirements. 

From a Stackelberg game perspective, the results reveal that the leader’s decisions can induce operating conditions that are not necessarily efficient from the follower’s dispatch perspective. In particular, during periods in which the low-cost resource is scarce, the exchange decisions determined at the upper level force the MV-DN operator to adopt a resource allocation that deviates from its independently optimal dispatch. As a consequence, the follower may face higher operating costs and less efficient outcomes in terms of generation allocation, voltage management, and congestion mitigation, which constitute its primary operational responsibilities. These results highlight the value of the proposed bilevel formulation, as it explicitly captures the strategic interaction between the leader and the follower and the resulting trade-offs between local and system-level objectives. Such behavior would not be adequately represented under the HPR, where the leader’s solution could still be optimal from its own perspective, but the operational inefficiencies imposed on the MV level would remain hidden. Therefore, the proposed bilevel approach provides a more comprehensive representation of the coordination problem, capturing both the interdependence of decisions across voltage levels and the potentially conflicting interests of the participating agents.

Although the proposed bottom-up bilevel framework captures hierarchical coordination between LV energy communities and the MV network while preserving discrete leader-level decisions, several extensions remain for future research. In the current formulation, discrete variables are limited to BESS operation, as the main objective was to validate the coordination structure and the LDD-ADMM algorithm under mixed-integer conditions. However, the framework can be naturally extended to include additional sources of discrete operational complexity at the LV level. In particular, future studies may incorporate P2P energy trading mechanisms, EV charging scheduling, demand-side flexibility programs, or other forms of local market participation. These additions would further increase the dimensionality and combinatorial nature of the upper-level problem, providing a richer representation of active distribution systems.

Finally, uncertainty was not explicitly modeled in the present study. The deterministic setting adopted here facilitates structural validation of the coordination scheme and its computational properties. Nonetheless, future work may incorporate stochastic or robust formulations at one or both levels. For example, uncertainty in PV generation forecasts, demand profiles, or wholesale market prices could be introduced at the MV level, while user-level flexibility uncertainty could be modeled within LV communities. Such extensions would enable evaluation of the resilience and robustness of the proposed coordination framework under realistic operating conditions. Additionally, the MV system was modeled as a closed operational layer, with no explicit interaction with the upstream transmission network. Future developments may integrate transmission-level constraints, locational marginal pricing signals, or participation in the MV-level reserve market. This would enable the study of cross-level coordination not only between LV and MV networks but also between distribution and transmission systems within a fully integrated hierarchical structure.


\section*{Acknowledgment}
This work has been supported by ANID FONDECYT Iniciación 11240745.

\appendix
\section{Single-Level Reformulation via KKT}\label{Appendix_A}
The reformulation of a bilevel optimization problem with a Stackelberg leader--follower structure into a single-level program via the optimality conditions of the lower-level problem requires that certain structural conditions be satisfied \cite{7942105}. In particular, the lower-level problem must admit an optimal solution for any feasible upper-level decision, be convex, and satisfy suitable regularity conditions ensuring the existence of dual variables and the validity of the KKT conditions as necessary and sufficient optimality conditions \cite{dempe2002foundations}.

In the proposed model, the feasibility of the lower-level MV-DN problem is guaranteed by construction. Infeasibility could only arise in two situations: if the MV-DN were unable to supply the total demand requested by the LV communities and the native MV demand, or if surplus energy injected by the LV communities could not be absorbed by any consumption node. The first situation is avoided by including sufficient dispatchable conventional generation at the MV level to meet the total system demand under all operating conditions. The second situation is prevented by the presence of additional MV demand nodes, which ensure that any surplus energy exported by the LV communities can always be assigned either to traditional MV loads or indirectly to other LV communities through the coordination mechanism.

Second, unlike top-down coordination schemes, all discrete decisions associated with battery operation are confined to the upper level. As a result, the MV-DN problem remains a continuous convex optimization problem, characterized by a convex quadratic objective function and linear equality and inequality constraints, with continuous decision variables only. Therefore, the lower-level problem is convex and satisfies strong duality. Under these conditions, the KKT optimality conditions characterize the lower-level solution. Consequently, the bilevel problem can be transformed into a single-level (SL) problem by embedding the primal feasibility, dual feasibility, and complementary slackness conditions of the MV-DN problem into the upper-level \cite{colson2007overview} formulation as follows:

\begin{equation}\label{SL_KKT}
\text{(SL)} \quad
\left\{
\begin{aligned}
\min_{\substack{\{x_k,\Delta b_k^L\}_{k \in \mathcal{K}},\\ y,\Delta b^M,\pi,\mu,\lambda}} \quad
&\sum_{k \in \mathcal{K}} F_k(x_k)
+ \sum_{k \in \mathcal{K}}\sum_{i \in \Omega_k^B}\sum_{t\in \mathcal{T}}
\left(
\lambda_t^c Mp^{bc}_{i,t}
+ \lambda_t^e Mp^{be}_{i,t}
- \lambda_t^s Lp^{s}_{i,t}
\right) \\[0.8em]
\text{s.t.} \quad
& G_k(x_k) = 0, \qquad \\
& H_k(x_k) \le 0, \qquad  \\
& \Delta b^L_{i,t}
= Mp^{bc}_{i,t} + Mp^{be}_{i,t} - Lp^{s}_{i,t},
\quad \forall i \in \Omega_k^B,\ \forall t \in \mathcal{T}, \\
& g(y) = 0, \\ 
& h(y) \le 0, \\ 
& \Delta b^M_{i,t} = -\Delta b^L_{i,t},
\quad \forall k \in \mathcal{K}, \forall i \in \Theta^k,\forall t \in \mathcal{T}, \\
& c + A^\top \pi + D^\top \mu + B^\top \lambda = 0, \\
& \mu \ge 0, \\
& \mu^\top h(y) = 0 .
\end{aligned}
\right.
\end{equation}

The SL formulation in \eqref{SL_KKT} replaces the lower-level optimization problem in \eqref{BL_problem} by its KKT optimality conditions. In particular, the primal feasibility of the MV-DN problem is preserved through constraints $g(y)=0$, $h(y)\leq 0$, and the coupling conditions $\Delta b^M_{i,t}= -\Delta b^L_{i,t}$, where the boundary exchanges are now explicit decision variables rather than fixed parameters. Dual feasibility is enforced by the non-negativity of the multipliers $\mu$, while complementary slackness is imposed by $\mu^\top h(y)=0$. Since the lower-level problem is a convex quadratic program with linear constraints, the stationarity condition can be written in compact matrix form as $c + A^\top \pi + D^\top \mu + B^\top \lambda = 0$, which corresponds to the evaluated gradient of the Lagrangian with respect to the MV-DN decision variables. Together, these conditions characterize the optimal response of the MV-DN and yield an exact single-level reformulation of the original bilevel problem.

\newpage
\section{Derivation of the LDD--ADMM coordination scheme}\label{Appendix_B}
This appendix details the derivation of the distributed coordination scheme from the bilevel formulation in \eqref{BL_problem}, leading to the LV and MV subproblems in \eqref{LV_ADMM_problem} and \eqref{MV_ADMM_problem}. Thus, to enable a distributed implementation, the coordination structure is reformulated by introducing duplicated variables associated with the boundary exchanges between the LVs and MV levels.
\begin{equation}\label{consistency_full}
Mp^{bc}_{k,t} = Lp^{bc}_{i,t}, \quad
Mp^{be}_{k,t} = Lp^{be}_{i,t}, \quad
Mp^{s}_{k,t} = Lp^{s}_{i,t}, 
\qquad \forall k \in \mathcal{K}, \forall i \in \Omega_k^B, \forall t \in \mathcal{T}.
\end{equation}


This reformulation separates local decision variables from shared coupling variables, allowing each LV network and the MV system to retain their internal information while interacting only through a reduced set of boundary quantities. As a result, the original bilevel problem \eqref{BL_problem} is recast into a structure more suitable for distributed coordination, referred to as the Decomposed Coordination of the Bilevel (DC-BL) problem, as defined below:

\begin{equation}\label{BL_reformulated_full}\small 
\text{(DC-BL)} \quad
\left\{
\begin{aligned}
\min_{\{x_k,\Delta b_k^L,Lp_k^{(\cdot)}\}_{k\in\mathcal{K}}} \quad
& \sum_{k\in\mathcal{K}} F_k(x_k)
+ \sum_{k\in\mathcal{K}}\sum_{i\in\Omega_k^B}\sum_{t\in \mathcal{T}}
\left(
\lambda_t^c Lp^{bc}_{i,t}
+ \lambda_t^e Lp^{be}_{i,t}
- \lambda_t^s Lp^s_{i,t}
\right) \\[0.6em]
\text{s.t.} \quad
& G_k(x_k)=0, \qquad \forall k\in\mathcal{K},\\
& H_k(x_k)\le 0, \qquad \forall k\in\mathcal{K},\\
& \Delta b^L_{i,t}
= Lp^{bc}_{i,t}+Lp^{be}_{i,t}-Lp^s_{i,t},
\qquad \forall i\in\Omega_k^B,\ \forall t\in \mathcal{T},\\[0.6em]
& (y,\Delta b^M, Mp^{(\cdot)}) \in
\arg\min_{y,\Delta b^M, Mp^{(\cdot)}}
\Big\{
f(y):\ g(y)=0,\ h(y)\le 0,\\
& \qquad\qquad\qquad\qquad
\Delta b^M_{k,t}
= -Mp^{bc}_{k,t}-Mp^{be}_{k,t}+Mp^{s}_{k,t},
\qquad \forall k\in\mathcal{K},\ \forall t\in \mathcal{T}
\Big\},\\[0.6em]
\multicolumn{2}{l}{\text{(Coupling constraints)}}\\[0.2em]
& Mp^{bc}_{k,t}=Lp^{bc}_{i,t},
\qquad \forall k\in\mathcal{K},\ \forall i\in\Omega_k^B,\ \forall t\in \mathcal{T},\\
& Mp^{be}_{k,t}=Lp^{be}_{i,t},
\qquad \forall k\in\mathcal{K},\ \forall i\in\Omega_k^B,\ \forall t\in \mathcal{T},\\
& Mp^{s}_{k,t}=Lp^{s}_{i,t},
\qquad \forall k\in\mathcal{K},\ \forall i\in\Omega_k^B,\ \forall t\in \mathcal{T}.
\end{aligned}
\right.
\end{equation}

The introduction of duplicated boundary variables assigns local copies to each LV network and the MV system, rendering their decision variables separable. As a result, the LV and MV problems can be solved independently, interacting only through the consistency constraints that link the duplicated variables. These equality constraints do not belong to either the LV or MV subproblems individually; rather, they define the coordination layer of the reformulated problem, ensuring consistency of the boundary exchanges without altering the underlying bilevel structure. In this sense, the original bilevel problem is reformulated as a distributed coordination problem over the boundary variables in \eqref{BL_reformulated_full}.

However, due to the definition of the net boundary exchange, $\Delta b^{L}_{i,t}
= Lp^{bc}_{i,t} + Lp^{be}_{i,t} - Lp^{s}_{i,t}$, and considering that the balance relation $\Delta b^{M}_{k,t} = -\Delta b^{L}_{i,t}$ is enforced within the follower problem, the three components are linearly dependent. Therefore, enforcing all consistency constraints in \eqref{consistency_full} is redundant, and one of them can be removed without altering the feasible set. In what follows, the coordination is imposed explicitly on the cheap-import and surplus-export components, while the expensive-import component is implicitly recovered from the balance relation. Under this reduced representation, the reformulated coordination problem can be interpreted as a structured consensus problem in which each LV network $k \in \mathcal{K}$ defines local copies of the boundary variables, while the MV system determines their global counterparts. 

In this context, the ADMM provides an effective framework to enforce consensus among distributed agents. In its standard form, a consensus problem can be written as
\begin{equation}\label{ADMM_standar}
\min \sum_{k \in \mathcal{K}} f_k(z_k) + g(z)
\quad \text{s.t.} \quad z_k = z,
\end{equation}
where $z_k$ and $z$ denote local and global copies of the shared decision variables, respectively. ADMM addresses this problem by constructing an augmented Lagrangian and decomposing it into a set of independent subproblems, corresponding in this case to each LV community and the MV system, which can be solved separately while being coordinated iteratively through the linking constraints $z_k = z$. In this way, consensus is enforced without requiring centralized information, as only the linking variables are exchanged among the agents.

To derive an ADMM-based coordination scheme from \eqref{ADMM_standar} for the reformulated problem \eqref{BL_reformulated_full}, an augmented Lagrangian must be constructed based on the identified coupling structure. It is important to emphasize that the augmented Lagrangian does not constitute an exact single-level reformulation of the bilevel problem in \eqref{BL_reformulated_full}. Instead, it is constructed from the reformulated coupling structure to enable a distributed coordination scheme that preserves the \textit{hierarchical interaction} between LV and MV levels. Accordingly, the augmented Lagrangian of the coordination problem is defined as
\begin{align}\label{augmented_lagrangian}
\mathcal{L} \;=\;& 
\sum_{k \in \mathcal{K}} F_k(x_k)
+ f(y)
+ \sum_{k \in \mathcal{K}} \sum_{t \in \mathcal{T}}
\mu^s_{k,t} \left( Mp^s_{k,t} - Lp^s_{k,t} \right)+ \sum_{k \in \mathcal{K}} \sum_{t \in \mathcal{T}}
\frac{\rho}{2}
\left\|
Lp^{bc}_{k,t} - Mp^{bc}_{k,t} + u^{bc}_{k,t}
\right\|^2,
\end{align}
It is worth noting that the augmented Lagrangian in \eqref{augmented_lagrangian} includes a quadratic penalty term only for the cheap-energy component, while the surplus component is handled through standard Lagrangian multipliers. Although both consistency constraints could be enforced via an augmented formulation, this would introduce additional consensus variables and increase the algorithm's computational burden. In this work, the surplus exchange is treated differently due to its structural role in the system. In particular, surplus energy injected by the LV networks must always be absorbed and redistributed by the MV system to ensure feasibility. As a result, this component does not require a strong consensus among agents and can be effectively coordinated via Lagrangian relaxation. Accordingly, the augmented Lagrangian in \eqref{augmented_lagrangian} is composed of the local objective functions of the LV networks, $F_k(x_k)$, and the objective function of the MV system, $f(y)$, together with the dual and penalty terms associated with the coordination constraints. The multipliers $\mu^s_{k,t}$ correspond to the Lagrangian relaxation of the surplus-exchange, while the quadratic penalty term, weighted by $\rho > 0$ and involving the scaled dual variables $u^{bc}_{k,t}$, enforces consensus on the cheap-energy component. 

The distributed coordination scheme is derived by exploiting the separable structure of the augmented Lagrangian in \eqref{augmented_lagrangian}. In particular, $\mathcal{L}$ is separable with respect to the LV and MV decision variables, except for the coupling terms associated with the shared boundary variables. This structure enables the application of a block-coordinate (Gauss--Seidel) minimization scheme, in which the augmented Lagrangian is minimized alternately with respect to each block of variables while keeping the remaining variables fixed at their most recent iterates \cite{7889039, li2026convergence}. Accordingly, for each iteration $r$, each LV subproblem is obtained by minimizing $\mathcal{L}$ with respect to the local variables $(x_k, b_k^L)$, while enforcing its original local constraints, i.e.,
\begin{equation}
\text{(LV-DN$_k$--ADMM)} \quad
\left\{
\begin{aligned}
\min_{x_k,\, b_k^L} \quad
& F_k(x_k)
+ \sum_{t \in \mathcal{T}} \sum_{i \in \Omega_k^B}
\left(
\lambda_t^{c} Lp^{bc}_{i,t}
+ \lambda_t^e Lp^{be}_{i,t}
- \lambda_t^s Lp^{s}_{i,t}
\right) \\
& \quad
+ \sum_{t \in \mathcal{T}} \sum_{i \in \Omega_k^B}
\frac{\rho}{2}
\left\|
Lp^{bc}_{i,t}
- Mp^{bc,(r)}_{k,t}
+ u^{bc,(r)}_{k,t}
\right\|^2 \\[0.5em]
\text{s.t.} \quad
& G_k(x_k) = 0, \\
& H_k(x_k) \le 0, \\
& \Delta b^{L}_{i,t}
= Lp^{bc}_{i,t}
+ Lp^{be}_{i,t}
- Lp^{s}_{i,t},
\quad
\forall i \in \Omega_k^B,\ \forall t \in\mathcal{T}.
\end{aligned}
\right.
\end{equation}

Similarly, the MV subproblem is obtained by minimizing $\mathcal{L}$ with respect to the variables $(y, \Delta b^M)$, subject to the original MV constraints.

\begin{equation}
\text{(MV-DN--ADMM)} \quad
\left\{
\begin{aligned}
\min_{y,\, \Delta b^M} \quad
& f(y)
+ \sum_{k \in \mathcal{K}} \sum_{t \in \mathcal{T}}
\mu^{s}_{k,t}
\left(
Mp^{s}_{k,t} - Lp^{s,(r)}_{k,t}
\right) \\
& \quad
+ \sum_{k \in \mathcal{K}} \sum_{t \in \mathcal{T}}
\frac{\rho}{2}
\left\|
Mp^{bc}_{k,t}
- Lp^{bc,(r)}_{k,t}
- u^{bc,(r)}_{k,t}
\right\|^{2} \\[0.5em]
\text{s.t.} \quad
& g(y) = 0, \\
& h(y) \le 0, \\
& \Delta b^{M}_{k,t}
= -\Delta b^{L,(r)}_{k,t},
\quad
\forall k \in \mathcal{K},\ \forall t \in\mathcal{T}.
\end{aligned}
\right.
\end{equation}

It is worth noting that the objective functions of the LV and MV subproblems are not structurally identical. Although both are derived from the same augmented Lagrangian, the resulting expressions differ due to the block-coordinate minimization scheme. In particular, for the cheap-energy variables coordinated via ADMM, the Gauss-Seidel decomposition implies that each block minimizes the same augmented Lagrangian with respect to its own variables while keeping the counterpart variables fixed at their most recent iterates. That is, for each LV$_k$ subproblem,
\[
\frac{\rho}{2} \left\| Lp^{bc}_{i,t} - Mp^{bc,(r)}_{k,t} + u^{bc,(r)}_{k,t} \right\|^2,
\]
where $Mp^{bc,(r)}_{k,t}$ is fixed, whereas in the MV subproblem,
\[
\frac{\rho}{2} \left\| Mp^{bc}_{k,t} - Lp^{bc,(r)}_{k,t} - u^{bc,(r)}_{k,t} \right\|^{2},
\]
where $Lp^{bc,(r)}_{k,t}$ is fixed. Hence, this corresponds to the same quadratic term viewed from each block with the other variable held constant, which explains why it appears in each subproblem. In contrast, the energy surplus is not handled via ADMM but rather via standard Lagrangian relaxation. As a result, its contribution remains linear and appears only in the MV subproblem. This is because, under the adopted decomposition, the corresponding dualized coupling constraint is enforced at the MV level, where the consistency of the surplus exchange is evaluated, and the associated dual variable is updated. The LV subproblems, in turn, retain the surplus variables as local decision variables but omit the corresponding dual term from their objective functions, resulting in an asymmetric representation of this term across subproblems. This modeling choice is further supported by the structural role of the surplus exchange in the system. Since the MV network is assumed to always absorb the surplus energy from the LV communities, the associated coupling does not require a strong consensus mechanism; a standard Lagrangian relaxation is sufficient to coordinate these exchanges through the corresponding dual signals.

\begin{table*}[!t]
\footnotesize
\begin{framed}
\begin{longtable}{p{3.5cm}p{10cm}}

\multicolumn{2}{l}{\textbf{Notation convention}}\\[2pt]
$(\cdot)^L$ & Element associated with the LV-DN.\\
$(\cdot)^M$ & Element associated with the MV-DN.\\[6pt]

\multicolumn{2}{l}{\textbf{Sets and indices}}\\[2pt]
$k \in \mathcal{K}$ & Set of LVs.\\
$t \in \mathcal{T}$ & Set of time periods.\\
$i \in \Omega_k$ & Set of buses in LV $k$.\\
$i \in \Omega_k^{A}$ & Subset of users at LV $k$; $\Omega_k^{A} \subseteq \Omega_k$.\\
$i \in \Omega_k^{B}$ & Subset of PCC buses at LV $k$; $\Omega_k^{B} \subseteq \Omega_k$.\\
$i \in \Theta$ & Set of buses in the MV network.\\
$i \in \Theta^{k}$ & Subset of MV buses with connected LV; $\Theta^{k} \subseteq \Theta$.\\
$i \in \Theta^{A}$ & Subset of MV buses with conventional demand and/or utility-scale resources; $\Theta^{A} \subseteq \Theta$.\\
$(i,j) \in \mathcal{L}_k$ & Set of lines in LV $k$; $i,j \in \Omega_k$.\\
$(i,j) \in \mathcal{L}$ & Set of lines in the MV network; $i,j \in \Theta$.\\[6pt]

\multicolumn{2}{l}{\textbf{Parameters}}\\[2pt]
$\lambda_t^{c}$ & Price of cheap energy imported by an LV-DN from the MV-DN.\\
$\lambda_t^{e}$ & Price of expensive energy imported by an LV-DN from the MV-DN.\\
$\lambda_t^{s}$ & Price of surplus energy exported by an LV-DN to the MV-DN.\\
$C^{pv}$ & Unit operating cost of PV generation.\\
$C^{bt}$ & Unit degradation cost of BESS discharging.\\
$R_{i,j}$ & Resistance of line $(i,j)$.\\
$X_{i,j}$ & Reactance of line $(i,j)$.\\
$PL_{i,t}$ & Active demand at bus $i$ and time $t$.\\
$QL_{i,t}$ & Reactive demand at bus $i$ and time $t$.\\
$Q^{(\cdot)}$ & Minimum and maximum reactive power limits.\\
$V_i^{(\cdot)}$ & Minimum and maximum squared voltage magnitude at bus $i$.\\
$I_{i,j}^{\max}$ & Maximum current magnitude in line $(i,j)$.\\
$\Gamma_i^{pv}$ & Installed PV capacity at bus $i$.\\
$\Gamma_i^{bt}$ & Installed BESS energy capacity at bus $i$.\\
$PV_t$ & Normalized solar irradiance profile at the LV level.\\
$PB_i$ & Maximum BESS charging/discharging power at LV bus $i$.\\
$SOC^{(\cdot)}$ & Minimum and maximum BESS state-of-charge limits.\\
$\phi^{(\cdot)}$ & BESS charging and discharging efficiencies.\\
$\nu_i$ & Binary parameter indicating whether a BESS is installed at bus $i$.\\
$PG_i$ & Maximum conventional generation capacity at MV bus $i$.\\
$C_i^{(\cdot)}$ & Quadratic ($a$), linear ($b$), and fixed ($c$) generation cost at bus $i$.\\[6pt]

\multicolumn{2}{l}{\textbf{Decision variables}}\\[2pt]
$p_{i,j,t}$ & Active power flow in line $(i,j)$ at time $t$.\\
$q_{i,j,t}$ & Reactive power flow in line $(i,j)$ at time $t$.\\
$\ell_{i,j,t}$ & Squared current magnitude in line $(i,j)$ at time $t$.\\
$v_{i,t}$ & Squared voltage magnitude at bus $i$ and time $t$.\\
$qg_{i,t}$ & Reactive power injection at bus $i$ and time $t$.\\
$\Delta p_{i,t}$ & Net active power injection at user bus $i$ and time $t$.\\
$\Delta b_{i,t}$ & Net boundary active power exchange at PCC bus $i$ and time $t$.\\
$pv_{i,t}$ & PV active power generation at bus $i$ and time $t$.\\
$pg_{i,t}$ & Conventional active power generation at bus $i$ and time $t$.\\
$ch_{i,t}$ & BESS charging power at bus $i$ and time $t$.\\
$ds_{i,t}$ & BESS discharging power at bus $i$ and time $t$.\\
$soc_{i,t}$ & BESS state of charge at bus $i$ and time $t$.\\
$w_{i,t}$ & Binary variable indicating BESS charging/discharging status.\\
$Mp_{k,t}^{s}$ & Power sold from the MV network to LV $k$.\\
$Mp_{k,t}^{be}$ & Expensive energy bought by LV $k$ from the MV network.\\
$Mp_{k,t}^{bc}$ & Cheap energy bought by LV $k$ from the MV network.\\

\end{longtable}
\end{framed}
\caption{Nomenclature.}\label{Nomenclature_Table}
\end{table*}
\newpage
\bibliographystyle{elsarticle-num-names}
\bibliography{references}

\end{document}